\documentclass[11pt]{article}

\usepackage{tikz-cd}
\usepackage{grffile}
\usepackage{ctable}
\usepackage{appendix}
\usepackage{geometry}
\usepackage{amssymb, nccmath, comment}
\usetikzlibrary{matrix,calc}
\usepackage{amsmath}
\usepackage{mathrsfs}
\usepackage[dvips]{epsfig}
\usepackage[small]{caption}
\usepackage{graphicx}
\usepackage[all]{xy}

\usepackage [english]{babel} 
\usepackage{enumerate}
\usepackage[utf8]{inputenc}
\usepackage[all]{xy}

\usepackage [autostyle, english = american]{csquotes}
\usepackage{graphicx}
\usepackage{bm}
\graphicspath{ {images4/} }
\MakeOuterQuote{"}
\newcommand{\Z}{\mathbb{Z}}
\newcommand{\N}{\mathbb{N}}

\newcommand{\C}{\mathbb{C}}
\newcommand{\R}{\mathbb{R}}

\newcommand{\ep}{\varepsilon}

\newcommand{\ds}{\displaystyle}

\newcommand{\T}{\mathcal{T}}
\newcommand{\V}{\Vert}
\newcommand{\lv}{\left\Vert}
\newcommand{\rv}{\right\Vert}

\newcommand{\wt}{\widetilde}
\newcommand{\wh}{\widehat}

\newcommand{\K}{\mathcal{K}}
\newcommand{\Id}{\mathbb{I}_n}
\newcommand{\coker}{\mathrm{coker}}
\newcommand{\lip}{\mathrm{Lip}}

\DeclareFontFamily{OT1}{pzc}{}
\DeclareFontShape{OT1}{pzc}{m}{it}{<-> s * [1.10] pzcmi7t}{}
\DeclareMathAlphabet{\mathpzc}{OT1}{pzc}{m}{it}

\newtheorem{Lemma}{Lemma}[section]
\newtheorem{Theorem}{Theorem}
\newtheorem{Proposition}[Lemma]{Proposition}
\newtheorem{Corollary}[Lemma]{Corollary}
\newtheorem{Remark}[Lemma]{Remark}

\newtheorem{Hypothesis}[Lemma]{Hypothesis}

\newenvironment{Proof}[1][\unskip]%
 {\begin{trivlist} \item[]{\bf Proof #1. }}%
 {\hspace*{\fill}$\rule{.4\baselineskip}{.4\baselineskip}$\end{trivlist}}

 {\begin{trivlist}\item[]\textbf{Acknowledgments.}}{\end{trivlist}}

 {\begin{center}\textbf{Abstract}}{\end{center}}

\makeatletter\@addtoreset{figure}{section}\makeatother

\makeatletter \@addtoreset{equation}{section} \makeatother

\makeatletter
\newsavebox{\@brx}
\newcommand{\llangle}[1][]{\savebox{\@brx}{\(\m@th{#1\langle}\)}%
  \mathopen{\copy\@brx\kern-0.5\wd\@brx\usebox{\@brx}}}
\newcommand{\rrangle}[1][]{\savebox{\@brx}{\(\m@th{#1\rangle}\)}%
  \mathclose{\copy\@brx\kern-0.5\wd\@brx\usebox{\@brx}}}
\makeatother

\definecolor{Green}{rgb}{0.,0.4,0.}

\renewcommand{\leq}{\leqslant}

\newcommand{\A}{\mathcal{A}}

\newcommand{\Rmnum}[1]{\uppercase\expandafter{\romannumeral #1\relax}}

\def\XXint#1#2#3{{\setbox0=\hbox{$#1{#2#3}{\int}$}
     \vcenter{\hbox{$#2#3$}}\kern-.5\wd0}}

\newfam\bifam
\font\tenbi=cmmib10 scaled \magstep1 \font\sevenbi=cmmib10 at 11pt
\font\fivebi=cmmib10 at 6pt \textfont\bifam = \tenbi
\scriptfont\bifam = \sevenbi \scriptscriptfont\bifam= \fivebi

\ifx\pdfoutput\undefined
   \pdffalse
\else
   \pdfoutput=1
   \pdftrue
\fi
\ifpdf
   \usepackage{graphicx}
   \usepackage{epstopdf}
\else
   \usepackage{graphicx}
\fi
\begin{document}

\begin{center}
{\fontsize{15}{15}\fontfamily{cmr}\fontseries{b}\selectfont{Coherent structures in nonlocal systems --- functional analytic tools}}\\[0.2in]
Olivia Cannon and Arnd Scheel\footnote{The authors acknowledge partial support through  grant NSF DMS-1907391.}\\[0.1in]
\textit{\footnotesize 
University of Minnesota, School of Mathematics,   206 Church St. S.E., Minneapolis, MN 55455, USA}
\date{\small \today} 
\end{center}

\begin{abstract}
We develop tools for the analysis of fronts, pulses, and wave trains in spatially extended systems with nonlocal coupling. We first determine Fredholm properties of linear operators, thereby identifying pointwise invertibility of the principal part together with invertibility at spatial infinity as necessary and sufficient conditions. We then build on the Fredholm theory to construct center manifolds for nonlocal spatial dynamics under optimal regularity assumptions, with reduced vector fields and phase space identified a posteriori through the shift on bounded solutions. As an application, we establish uniqueness of small periodic wave trains in a Lyapunov center theorem using only $C^1$-regularity of the nonlinearity.
\end{abstract}

\begin{center}
    \textbf{Keywords:} Nonlocal, center manifold, Fredholm, coherent structures, Lyapunov center
\end{center}

\setlength{\parskip}{4pt}
\setlength{\parindent}{0pt}

\section{Introduction}

Describing the emergence of coherent structures and self-organized collective behavior in large complex systems is both central to our understanding of dynamical behavior and theoretically challenging. Recently, interest has grown in studying systems with nonlocal coupling, with motivation from neuronal networks, biology, material science, and ecology  \cite{NCalfimov1992dynamical,NCberner2020solitary,NCsilva2014characterization,NCdu2010nonlocal, NCfaye2015modulated,NCgopal2014observation,NCjaramillo2022rotating,NCpalmer2001nonlinear,NCtanaka2003complex,NCtang2018nonlinear,truong2021global,NCyonker2006nonlocal,zhou2006materials}. Nonlocal coupling can take many forms but, in a continuum modeling context, can be quite generally represented by integral operators, rather than differential operators in local differential equation models. Phenomena in nonlocally coupled systems are often qualitatively different from phenomena in differential  equations, notably including effects ranging from 
singularity formation \cite{hur2012singularities,seleson2009peridynamics,zhou2006materials}, to rapid synchronization \cite{chiba2018kuramoto}, pinning \cite{anderson2016pinning},  or acceleration of fronts \cite{cabre2009fraclaplacian,meleard2015singular}. Relatedly, mathematical techniques from differential equations are not immediately applicable to nonlocally coupled systems and limitations of techniques often point to new phenomena \cite{anderson2016pinning,chiba2015proof,truong2021global}.


In the present work, we focus on developing techniques that adapt tools from the study of differential equations to nonlocal systems, identifying in particular potential limitations such as the lack of regularity or the loss of compactness. Our focus is on coherent structures, particularly traveling waves---periodic wave trains, solitary waves and pulses,  fronts, and other types of solutions arising from the inherent self-organizing capabilities of large systems. In the analysis of existence, stability, and bifurcations of such states, one desires a robust functional-analytic framework which the present work aims to contribute to. The class of equations we study arises as steady-state or traveling-wave equations of a time-dependent system. Our contributions can be organized into three categories. We study, for a class of nonlocal equations:
\begin{itemize}
    \item Fredholm theory for linearization at coherent structures;
    \item Center manifold theory for bifurcation of coherent structures from the trivial state;
    \item A Lyapunov-Center theorem for nonlocal systems.
\end{itemize}
In fact, technical results in these three areas build on each other, with center manifold theory relying on Fredholm theory, and uniqueness in Lyapunov-Center theorems relying on center manifold theory. We describe the above contributions in more detail below, briefly summarizing results as well as connections to local theory. 

\paragraph{Fredholm Properties.}

Fredholm theory is instrumental in the study of bifurcation theory in local as well as nonlocal settings. For instance, in situations where a given coherent structure exists but the linearization of the system is not invertible, a Fredholm linearization may allow one to establish continuation and bifurcation results using Lyapunov-Schmidt reduction. Fredholm properties for a nonlocal operator on $L^2(\R,\C^n)$ corresponding to traveling wave solutions of a time-dependent nonlocal equation were first established in \cite{faye2014fredholm}. Here, we consider instead the somewhat broader class of operators of the form 
\begin{equation}\label{e:operator} \T U(\xi) = A(\xi)U(\xi) + K_\xi*U,\end{equation}
with $\xi \in \R$, $U(\xi) \in \C^n$, corresponding to steady-state solutions of a time-dependent nonlocal equation. As in \cite{faye2014fredholm}, the convolution kernel is \textit{inhomogeneous}, $\xi-$dependent, with limits at $\pm\infty$, but possesses some smoothing properties. In contrast to  \cite{faye2014fredholm}, the principal part of the operator is a multiplication operator rather than a differential operator, and we investigate Fredholm properties of $\mathcal{T}$ on a larger class of function spaces, $L^p(\R,\C^n), 1 \le p \le \infty$, as well as $C^0(\R, \C^n)$, the latter with an eye toward proving center manifold properties in later sections. We note however that Fredholm properties are useful beyond the study of small-amplitude structures: they have been used to investigate eigenvalue problems near the edge of the essential spectrum \cite{faye2014fredholm} or to construct a Conley-Floer homology theory for gradient-like problems 
\cite{bakker2019large} and thus establish existence of large-amplitude front solutions for nonlocal systems.  

Our results identify necessary and sufficient conditions for the operators of the form  \eqref{e:operator} to be Fredholm (Theorem \ref{main}), and show how to compute the index (Theorem \ref{indexflow}). Informally, Theorem \ref{main} states that within a large class of operators, 
\[
\T \textrm{ is Fredholm } \iff \begin{cases} \T \textrm{ is invertible at spatial infinity;} \\ \textrm{The multiplication operator } A(\cdot) \textrm{ is invertible.}  \end{cases} 
\]
The first condition, loss of invertibility at infinity, is a well-known source of non-compactness, also in local problems \cite{faye2014fredholm,palmer1988fredholm}. The second condition arises from the change in principal part and loss of regularity, and may contribute to possible bifurcations such as depinning of fronts in the nonlocal setting \cite{anderson2016pinning} or synchronization transitions in coupled oscillators \cite{chiba2016mean}. 

\paragraph{Nonlocal Center Manifolds.}
Center manifold theory has long been used to study small-amplitude solutions of nonlinear equations. Originally set in finite dimensions \cite{Kelley1967TheSC}, then extended to Banach space settings \cite{Henry1989GeometricTO} and ill-posed equations \cite{kirchgassner1982wave}, the reduction of large or infinite-dimensional systems to a low-dimensional submanifold can allow for, for instance, existence and uniqueness arguments where they otherwise are not possible. We are concerned here with the construction of small, bounded stationary or traveling-wave solutions for nonlinear, nonlocal equations. For local equations, for instance PDEs set on $x\in\R$ or $x$ in a cylinder, such solutions can be studied using spatial dynamics and existing center manifold results. Constructing such stationary or traveling-wave solutions for nonlocal equations poses new challenges, in particular since an initial-value problem formulation, even an ill-posed one as for elliptic equations, is not readily available. Analytical results therefore were limited to special kernels that allow for a reformulation as an ODE \cite{faye2013existence}. This obstruction was removed in  \cite{fs2018center}, with a center manifold theory for nonlocal systems of the form 
\begin{equation}
0 = -U + K*U + F(U), 
\end{equation}
for exponentially localized $K$. There, the need for a phase space is sidestepped: instead of parameterizing initial conditions over a center subspace, entire \textit{trajectories} are parameterized  in function space over the kernel of the linearization, which is finite-dimensional. The crux of this idea is that the analogue of a flow in phase space is the shift operator in function space---the shift operator $\tau_\xi$ ``flows''
a trajectory $u(\cdot)$ forward to the shifted trajectory $u(\cdot + \xi)$. This flow, the action of the shift operator, can then be pulled back to the kernel and differentiated, in order to obtain a reduced vector field. 

Along with establishing a center manifold comes the question of optimal regularity. In traditional settings, one seeks to establish $C^k$ regularity of center manifolds for $C^k$ vector fields, or $C^{k,\alpha}$ regularity for $C^{k,\alpha}$ vector fields  \cite{Henry1989GeometricTO}, for finite-dimensional or Banach space settings \cite{bialy2002Ck,Gallay1993ACM,Kelley1967TheSC,sandstede2018Ck,VdBVH1987}. The phrase "$C^k$ manifold" refers, equivalently in that case, to regularity of the map parameterizing the set of center solutions, as well as to the regularity of the reduced vector field. Analogous 'optimal regularity' results were precluded in  \cite{fs2018center}, through the use of an $H^1$-function space setting: although the proof there establishes $C^k$ regularity of the manifold for a $C^k$ nonlinearity on $H^1$, a pointwise nonlinearity must be a $C^{k+1}$ function in order for the substitution operator to be a $C^k$ operator on $H^1$. We remedy this loss of regularity by relying on a  $C^0$ function-space setting, where pointwise substitution operators do not lose regularity.

Our contribution in Theorem \ref{CMexistence} then is twofold: 
\begin{itemize}
    \item \emph{Optimal regularity:} We construct nonlocal center manifolds on $C^0$ spaces, yielding $C^k$ manifolds and reduced vector fields for $C^k$ pointwise nonlinearities after a $C^k$ change of coordinates. 
    \item \emph{Local cutoff:} Our construction on $C^0$ spaces does not rely on the modified cutoff function necessary in the $H^1$-setting \cite{faye2020corrigendum}, simplifying the argument and allowing easier adaptation to different nonlinearities.
\end{itemize}
We delineate in Section 4 of this paper this construction of the nonlocal center manifold on $C^0$-based spaces, a key ingredient of which is the Fredholm theory from Sections 2 and 3. We note that the $C^k$ change of coordinates is not necessary to achieve the $C^k$ map parameterizing the center manifold, only the reduced vector field, since it allows bootstrapping of the center solutions. We remark also that we are able to recover the smoothness from \cite{fs2018center} when changing back to the original coordinates. Of course, the existence of a $C^k$ reduced vector field in \textit{some} coordinates, may well be useful and desirable since it allows for arguments based on uniqueness or sharp Taylor expansions that yield results which are valid independent of coordinates choices.

\paragraph{Lyapunov-Center Theorem.} In Hamiltonian and reversible systems, one can often conclude the existence of nonlinear oscillations from oscillations in the linear part. Such Lyapunov-Center theorems have been established in many contexts \cite{buzzi2004reversible,devaney76LC,li2011hamiltonian,schmidt1978center}.  First, existence of a one-parameter family of periodic trajectories near an equilibrium of the nonlinear flow is guaranteed by a pair of nonresonant imaginary eigenvalues $\pm i \omega$. Uniqueness of this family within the class of small periodic solutions can then be guaranteed by Lyapunov-Schmidt reduction, if there is exactly one simple pair of imaginary eigenvalues. Further, if a center manifold exists, one can show uniqueness of the family within the class of all small bounded, not necessarily periodic solutions to the nonlinear system. 

As an application of the center manifold on $C^0$ spaces, we prove here a Lyapunov-Center theorem for a system 
\begin{equation}\label{e:LC1} 0 = -u + k*(Au + N(u)) \end{equation}
with $A$ a constant matrix, $k$ an exponentially localized kernel, and $N(u)$ a $C^1$ pointwise nonlinearity, $N(0) = N'(0) = 0$. In the nonlocal case, reversibility corresponds to evenness of the convolution kernel $k$, and eigenvalues in the classical systems correspond to roots of the equation $0 = d(\nu) =\det(\Id + \wh{k}(\nu)A)$. Our result, Theorem \ref{LC}, can thus informally be stated:  
\[
\begin{cases} 
d(\nu) \textrm{ has a unique pair of roots }\pm i \omega_* \textrm{ on } i\R, \\ d'(i\omega_*) \neq 0 \end{cases} 
\implies \textrm{ all small bounded solutions to \eqref{e:LC1} are periodic.} 
\] 
In the context of spatial dynamics, the result establishes absence of small-amplitude coherent structures, such as solitary waves or nanopterons, for wave speeds different from group velocities under a non-resonance condition, with optimal regularity assumptions; see Remark \ref{r:LC}.

Technically, the $C^1$ case requires careful analysis because the principal term in the reduced equation is essentially quadratic. The argument also relies on the ability to establish a center manifold for a $C^1$ nonlinearity, which is made possible by the $C^0$-based center manifold construction.

\paragraph{Outline of the Paper.}
We establish in Section 2 necessary and sufficient conditions for Fredholm properties of a class of nonlocal operators. We characterize Fredholm indices of these operators in Section 3 via relative Morse indices, requiring stronger localization of the kernel than in the previous section. 

In Section 4, we prove existence of center manifolds for nonlocal systems on $C^0$-based spaces, using the methods in \cite{fs2018center}, and establish regularity of the reduced vector field. We use this center manifold theory in Section 5 to prove a Lyapunov-Center theorem for nonlocal systems. 

We note that sections may be read independently from each other---taking major results of the others for granted, and occasionally notation, they are essentially self-contained.

\section{Fredholm Properties of a Nonlocal Operator}\label{s:2}

We establish Fredholm properties for a class of nonlocal operators whose principal part is a multiplication operator, with a lower-order integral operator perturbation. Such operators may arise in studying space-dependent equilibria of nonlocal differential equations. The nonlocal coupling here is not a true convolution, except in the limit at spatial infinity. Such operators arise from linearizations at heteroclinic profiles in otherwise translation invariant systems, or from translation-invariant problems considered in weighted spaces. Our main focus in this section is to establish
Fredholm properties for these operators on $L^\infty$ and $C^0$, and outline adaptations to $L^p, p \ge 1$. In subsequent sections, Fredholm properties on $C^0$ will be used to extend results in \cite{fs2018center} on nonlocal center manifolds to $C^0$-based spaces. This class of operators are related to those in \cite{faye2014fredholm} but slightly more general and interesting in their own right due to additional sources of loss of compactness. 
 
\subsection{Setup and notation}\label{s:2.1}

We denote by $L^1, L^\infty$ the usual $L^p$ spaces $L^1(\R,\C^n)$ and $L^\infty(\R, \C^n)$, and we let $C^0(\R, \C^n)$ be the space of continuous functions with finite norm,
\[\lv f \rv_{C^0} = \max_{1\leq i\leq n}\sup_{x \in \R} | f_i(x)|.
\]
We let $M_n(\C)$ be the set of $n \times n$ complex matrices. Lastly, we introduce the weighted $L^2$ space \[
L^2_1(\R, M_n(\C))=
\{f \in L^2(\R, \C^n) \ | \ \| 
\sqrt{1 + \xi^2} \cdot f \|_{L^2} < \infty \}.\] Also define the complex Fourier Transform on $L^2(\R, \C^n)$ by
\[\wh{f}(i\ell) = \frac{1}{\sqrt{2\pi}}\int_\R f(\xi)e^{- i\ell \xi}d\xi;
\]
note that the standard Fourier transform evaluates $\wh{f}$ on $i\R$.
\paragraph{A class of nonlocal operators.}

We consider operators of the form 
\begin{equation}\label{e:opdef} \begin{split} \T: L^\infty(\R, \C^n) &\to L^\infty(\R,\C^n) \\
U(\xi) &\mapsto A(\xi)U(\xi) + \int_{\R}K(\xi - \xi'; \xi)U(\xi')d\xi',
\end{split}
\end{equation}
with $A(\cdot) \in L^\infty(\R, M_n(\C))$, $K(\xi -\cdot, \xi) \in W^{1,\infty}(\R, W^{1,1}(\R, M_n(\C))$. We further require, in order to ensure properties of the adjoint, that $K(\cdot + \xi, \cdot) \in W^{1,\infty}(\R, W^{1,1}(\R,M_n(\C))).$ The integral kernel $K_\xi(\cdot) = K(\cdot, \xi)$ can be thought of as an inhomogeneous convolution $K_\xi * U$. We denote the pair $(A,K)=: \A$, and we denote the operator by $\T_\A$ or simply $\T$, when unambiguous. 
We also consider the analogous class of operators $\T_{\A^C}$ on $C^0(\R, \C^n)$, for which all assumptions are the same except that we must have $A(\cdot) \in C^0(\R, M_n(\C))$. 

With suitable assumptions on limits and regularity, we wish to establish Fredholm properties for operators in this class: we first identify necessary and sufficient conditions for the operator to be Fredholm and then, with stronger localization assumptions, relate the Fredholm index to a spectral flow. 

We give two examples where such generalized convolution kernels arise.
\paragraph{Linearized nonlocal Allen--Cahn equation.}  The reaction-diffusion Allen--Cahn equation can be posed nonlocally as \begin{equation} \frac{du}{dt} = d(-u + k*u) + f(u), \end{equation} 
with $\int k=1$, for instance a normalized Gaussian, and constant effective diffusivity $d>0$. Considering stationary solutions, and then linearizing the equation about an interface-like solution $u_*$, $u_*(x)\to u_\pm$ for $x\to \pm\infty$, one obtains
\begin{equation} \T u = d(-u + k*u) + f'(u_*)\cdot u. \end{equation}

\paragraph{Neural Fields.}
Similarly, one can consider simple models for neural fields with an assumption of homogeneity, that is, translation invariance, 
\begin{equation} \frac{du}{dt} = -u + k*F(u), 
\end{equation} 
where $x$ lives in physical or feature space, and $u$ denotes a possibly averaged  state of the neural field. The state could  be scalar- or vector-valued, and the convolution kernel is often assumed Gaussian, or, for technical reasons,  to possess  rational Fourier transform. 

Again considering stationary solutions $u_*$ and linearizing, one obtains
\begin{equation}\T u =-u + k*[F'(u_*)\cdot u]. 
\end{equation} 


\subsection {Fredholm Properties}\label{s:2.2}
We state the main result and hypotheses.
\begin{Hypothesis}[Limits at Infinity]\label{hyp1}

We assume that there exist two matrices $A^\pm \in M_n(\C)$ such that \[A(\xi) \to A^\pm, \  \xi \to \pm \infty.\]

We also assume that there exist two functions $K^{\pm} \in W^{1,1}(\R, M_n(\C)) \cap L^2_1(\R, M_n(\C))$ such that \[ \lim_{\xi \to \pm \infty} \lv K(\cdot, \xi) - K^\pm(\cdot)\rv_{L^1} = 0 \] 
\[ 
\textrm{and } \lim_{\xi \to \pm \infty} \lv K( \cdot  ,  \cdot + \xi) - K^\pm(\cdot)\rv_{L^1} = 0. 
\] 
\end{Hypothesis}


\begin{Remark} The examples discussed above can be shown to satisfy these hypotheses, given somewhat mild assumptions on the kernel, and assuming  that $u_*\in L^\infty$, with limits at infinity, so that $f'(u_*(x)) \to f'(u^\pm(x)), x \to \pm \infty$. 
\end{Remark}

\begin{Theorem}\label{main} Given $
\T_\A$ in the class of operators defined in Section \ref{s:2.1} that satisfies Hypothesis \ref{hyp1},  the following are equivalent: 
\begin{enumerate}[(i)]
\item  $\T_\A$ is Fredholm.
\item $\T_\A$ satisfies
\begin{itemize}
    \item[(a)] {Hyperbolicity at Infinity:} $\det(\widehat{K^\pm}(i\ell) + A^\pm) \neq 0 $ for all $\ell\in\R$;
    \item[(b)]{Pointwise Invertibility of Principal Part:} $A^{-1}(\cdot)\in L^\infty(\R, M_n(\C))$.
\end{itemize}

\end{enumerate}

\it
Furthermore, when $\T_\A$ is Fredholm, its index depends only on the limits $A^\pm$ and $ K^\pm(\cdot)$ defined in Hypothesis \ref{hyp1}. 

If $A$ is continuous, the analogous result holds for $\T_{\A^C}$.

\end{Theorem}


Note that the above theorem essentially suggests that loss of compactness happens for two reasons---noncompactness of the domain, and pointwise lack of invertibility in the principal part. The former is a well-known source of non-compactness, but the latter is not present in the results of \cite{faye2014fredholm}, where the principal part is a differential operator. 

\subsection{Proof of Theorem 1, Sufficiency of Conditions for Fredholm} 

We first state some propositions which will be used in the proof.

\begin{Proposition}\label{bddH-1} There exists $C > 0$ such that for all $U \in L^\infty$, the following estimate holds: \begin{equation}\label{bddH} \left\Vert  \int_{\R}K(\xi - \xi'; \xi)U(\xi')d\xi' \right\Vert_{L^\infty} \le C \Vert U \Vert_{(W^{1,1})^*}\,.\end{equation}\end{Proposition}

\begin{Proof} Consider the operator \begin{align*}
\wt{\mathcal{K}}: L^1(\R, \C^n) &\to W^{1,1}(\R, \C^n) \\
U(\xi) &\mapsto \int_\R K^*(\xi' - \xi, \xi')U(\xi')d\xi',
\end{align*} where $K^*$ denotes the conjugate transpose of the matrix $K$. 

We have that \begin{align*}\lv \wt{\mathcal{K}}U\rv_{W^{1,1}} &\le \left( \sup_{\xi} \lv  K^*(\cdot + \xi, \cdot)\rv_{L^1} + \sup_{\zeta} \lv  \frac{d}{d\xi}K^*(\xi- \zeta, \xi)\rv_{L^1(\xi)}\right)\lv U \rv_{L^1}\\
&\le C\lv U \rv_{L^1},
\end{align*} by the fact that $K(\cdot + \xi, \cdot) \in W^{1,\infty}(\R, W^{1,1}(\R,M_n(\C))).$ Therefore, $\wt{\K}$ is a bounded operator from $L^1(\R, \C^n)$ to $W^{1,1}(\R, \C^n)$.
Since $\wt{\K}$ is bounded,  its adjoint $\wt{\K}^*$ is also bounded as an operator from $(W^{1,1}(\R, \C^n))^*$ to $L^\infty(\R, \C^n)$.

Formally, $\wt{\K}^*$ is defined only as an abstract operator on $(W^{1,1})^*$; however, elements of $L^\infty$ can be considered elements of $(W^{1,1})^*$ via the measure associated to the $L^\infty$ function. Whenever the argument of $\wt{\K}^*$ corresponds to an $L^\infty$ function in this way, the adjoint operator $\wt{K}^*$ must coincide with the operator 
\begin{equation}\label{e:formaladj} U(\xi) \to \int_{\R}K(\xi - \xi'; \xi)U(\xi')d\xi'.\end{equation}

The boundedness of the adjoint operator then implies the boundedness of the operator \eqref{e:formaladj} on $L^\infty$, which gives for $U \in L^\infty$ the bound in \ref{bddH-1} as desired. \end{Proof}

We will also need the following lemma:

\begin{Lemma}[Abstract Closed Range Lemma]\label{AbsClRg} Suppose that $X, Y,$ and $Z$ are Banach spaces, that $\T$ is a bounded linear operator, and that $\mathcal{R}: X \to Z$ is a compact linear operator. Assume that there exists a constant $c >0$ such that 
\[
\lv U \rv_X \le c\left(\lv \T U \rv_Y + \lv \mathcal{R} U \rv_Z\right), \ \ \ \  \textrm{ for all } U \in X.
\]
Then $\T$ has closed range and finite-dimensional kernel. 
\end{Lemma}

\begin{Proof} See \cite{schwarz1993morse}. \end{Proof}

\begin{Proposition}\label{mainprop} For the operator $\T_\A$, there exist constants $c > 0$ and $L >0$ so that \begin{equation}\label{prop} \lv U \rv_{L^\infty} \le c(\lv U\rv_{(W^{1,1}([-L,L], \C^n))^*} + \lv \T_\A U \rv_{L^\infty}).\end{equation} The same holds for $\T_{\A^C}.$\end{Proposition}
\rm In particular, this will allow us to use Lemma \ref{AbsClRg} since the composition $\mathcal{I} \circ \mathcal{R}_L$ of the restriction operator $\mathcal{R}_L: L^\infty(\R, \C^n) \to L^\infty([-L,L],\C^n)$ and the inclusion operator $ \mathcal{I}: L^\infty([-L,L],\C^n) \to (W^{1,1}([-L,L],\C^n))^*$ is compact. The same argument applies for the  analogous operators on $C^0$ since the latter is a closed subspace of $L^\infty$.

\begin{Proof}
Let $\T$ refer either to $\T_\A$ or $\T_{\A^C}$. Following \cite{faye2014fredholm}, we divide the proof into four steps. 

\it Step 1: \rm  We first show that \begin{equation} \label{step1} \lv U\rv_{L^\infty} \le c_1(\lv U\rv_{(W^{1,1}(\R, \C^n))^*} + \lv \T U\rv_{L^\infty}).\end{equation}

For each $U$, we have 
\begin{align*} \lv \T U \rv_{L^\infty} &= \left\Vert A(\xi)U(\xi)+ \int_{\R}K(\xi - \xi'; \xi)U(\xi')d\xi'\right\Vert_{L^\infty}\\ &= \lv A(\xi)\left(U(\xi) + A^{-1}(\xi)\int_\R K(\xi - \xi'; \xi)U(\xi')d\xi' \right) \rv_{L^\infty} \\ &\ge \frac{1}{\Vert A^{-1}(\xi)\Vert_{L^\infty}}  \left\Vert U(\xi) + A^{-1}(\xi)\int_{\R}K(\xi - \xi'; \xi)U(\xi')d\xi' \right\Vert_{L^\infty}\\ 
&\ge c\left\Vert U(\xi) + A^{-1}(\xi)\int_{\R}K(\xi - \xi'; \xi)U(\xi')d\xi' \right\Vert_{L^\infty}\\
&\ge c\left(\Vert U(\xi)\Vert_{L^\infty} - \V A^{-1}(\xi) \V_{L^\infty} \cdot \left\Vert \int_{\R}K(\xi - \xi'; \xi)U(\xi')d\xi' \right\Vert_{L^\infty}\right) \\
& \ge c(\Vert U(\xi)\Vert_{L^\infty} - c'\Vert U(\xi)\Vert_{(W^{1,1})^*}),
\end{align*}

for some $c, c' >0$, which implies the estimate \eqref{step1}.

\it Step 2: \rm We now consider a constant-coefficient operator \begin{equation}\label{coco} \T_{\A^0}: U(\xi) \mapsto A^0U(\xi) + (K^0*U)(\xi), \end{equation}
with $A^0$ invertible and $K^0(\cdot) \in W^{1,1}(\R, M_n(\C))  \cap   L^2_1(\R, M_n(\C))$, satisfying the hyperbolicity condition $\det(\wh{K^0}(i\ell) + A^0) \neq 0 $ for all $ \ell \in \R.$ Note that the condition $K \in L^2_1(\R, M_n(\C))$ guarantees that $\wh{K} \in H^1$. We will show directly that $\T_{\A^0}$ is bounded invertible, since $L^\infty$ is less amenable to the properties of Fourier multipliers. 

We define the inverse of $\T_{\A^0}$ on $ L^\infty(\R, M_n(\C))$ by \begin{equation}(\T_{\A^0})^{-1}(U) = (A^0)^{-1}U + K_{\textrm{inv}}*U,\end{equation} where $K_{\textrm{inv}}$ is the inverse Fourier transform of $ ((A^0+\wh{K^0}(i\ell))^{-1} - (A^0)^{-1})$.

That $\T_{\A^0}^{-1}$ is an inverse can be shown directly by calculating \[\T_{\A^0}(\T_{\A^0})^{-1}U = U + (A^0K_{\textrm{inv}} + K^0A^{-1} + (K*K_{\textrm{inv}}))*U;\] the function $(A^0K_{\textrm{inv}} + K^0A^{-1} + (K*K_{\textrm{inv}}))$ is the Fourier inverse of the 0 function, so vanishes almost everywhere, giving $\T_{\A^0}(\T_{\A^0})^{-1}U = U$. One can likewise show the same for $(\T_{\A^0})^{-1}\T_{\A^0}$. We would then like to show that the inverse function $(\T_{\A^0})^{-1}$ is bounded. We do this by showing that $\wh{K_{\textrm{inv}}} \in H^1$, so that $K_{\textrm{inv}} \in L^1$, leading to boundedness of $\T_{\A^0}^{-1}$.

We have by a matrix identity that $(A^0+\wh{K^0}(i\ell))^{-1} - (A^0)^{-1} = -(A^0)^{-1}\wh{K^0}(i\ell)(A^0+\wh{K^0}(i\ell))^{-1}$. Then \begin{align*}\lv((A^0+\wh{K^0}(i\ell))^{-1} - (A^0)^{-1})\rv_{M_n(\C)} &\le \left(\sup_{\ell \in \R}\lv(A^0)^{-1}\left(A^0+\wh{K^0}(i\ell)\right)^{-1}\rv_{M_n(\C)}\right)\lv \wh{K^0}(i\ell)\rv_{M_n(\C)}\\ 
&= C_1 \lv \wh{K^0}(i\ell) \rv_{M_n(\C)},\end{align*}  and 
\begin{align*}\lv \frac{d}{d\ell}((A^0+\wh{K^0}(i\ell))^{-1} - (A^0)^{-1})\rv_{M_n(\C)}&\le\sup_{\ell \in \R} \left( \lv \left(A^0+\wh{K^0}(i\ell)\right)^{-2} \rv_{M_n(\C)}\right)  \lv (\wh{K^0})'(i\ell) \rv_{M_n(\C)}\\
&= C_2 \lv (\wh{K^0})'(i\ell) \rv_{M_n(\C)},
\end{align*}
so that $\lv \wh{K_{\textrm{inv}}}\rv_{H^1} \le (C_1 + C_2) \lv \wh{K^0} \rv_{H^1}$, and $\wh{K_{\textrm{inv}}} \in H^1$. Then $K_{\textrm{inv}} \in L_1^2(\R, M_n(\C)) \subset L^1(\R, M_n(\C))$.

For $c_2 :=  \left(\lv (A^0)^{-1}\rv_{M_n(\C)} + \lv K_{\textrm{inv}} \rv_{L^1}\right)$, this gives the estimate

\begin{equation} \lv U \rv_{L^\infty} \le c_2\lv\T_{\A^0} U\rv_{L^\infty}, \ \ \ \textrm{ for all }  U \in L^\infty. \end{equation}

\it Step 3: \rm We now want to show that there exists $L >0$ such that if $U(\xi) = 0$ for $|\xi| < L-1, $ we have \begin{equation}\label{step3}\Vert U\Vert_{L^\infty} \le c_3\Vert \T U \Vert_{L^\infty}. \end{equation}

First, suppose that we have two functions $U^+(\xi) = 0, \xi \le L-1$, and $U^{-}(\xi) = 0, \xi > -(L-1)$. 

Then, note that since $K, A$ satisfy Hypothesis \ref{hyp1}, we may find $L$ large enough, so that for $U^\pm$, 
\[
\lv \int_\R\left(K(\xi - \xi';\xi) - K^\pm(\xi - \xi')\right)U^\pm(\xi')d\xi' \rv_{L^\infty} \le \frac{\ep}{2} \V U^\pm \V_{L^\infty} 
\]
\[ 
\V(A^\pm - A)U^\pm \V_{L^\infty} \le \frac{\ep}{2} \V U^\pm \V_{L^\infty} ,
\]
so we have $\frac{1}{c_2} \V U^\pm \V_{L^\infty} \le \V \T^\pm U^\pm \V_{L^\infty} \le \ep \V U^\pm \V_{L^\infty} + \V \T U^\pm \V_{L^\infty}$, which gives $\Vert U^\pm\Vert_{L^\infty} \le c \Vert \T U^\pm \Vert_{L^\infty}$, choosing $\ep c_2 < 1$, where the implicit notation $\T^\pm$ refers to the map 
\[
\T^\pm: U(\xi) \mapsto A^\pm(\xi)U(\xi) + \int_{\R}K^\pm(\xi - \xi')U(\xi')d\xi').
\]
Finally, given $U$ such that $U = 0, |\xi| < L-1$, we can decompose $U = U^+ + U^-$, taking $U^+(\xi) = U(\xi), $ for $ \xi > 0$, and $U^+(\xi) = 0$ otherwise, and $U^-(\xi) = U(\xi)$ for $\xi \le 0$, and $U^-(\xi) =$ 0 otherwise. 

Then we have \[\V U \V_{L^\infty} \le \V U^+ \V_{L^\infty} + \V U^- \V_{L^\infty} \le c(\V \T U^+ \V_{L^\infty} + \V \T U^- \V_{L^\infty}) \le 2 c \V \T U \V_{L^\infty} =: c_3 \V \T U \V_{L^\infty} ,\] as desired.

\it Step 4: \rm Finally, let $\chi$ be a smooth cutoff function equal to 0 outside $[-L,L]$ and equal to 1 for $|\xi| < L-1$. 
Then we have 
\begin{align*}
\V U \V_{L^\infty} &\le \V \chi U \V_{L^\infty} + \V (1 - \chi)U \V_{L^\infty} \\
&\le c_1(\V \chi U \V_{W^{1,1^*}} + \V \T(\chi U) \V_{L^\infty}) + c_3\V \T (1 - \chi)U \V_{L^\infty} \textrm{ \ \ \ (by steps 1 and 3 )}\\
&\le c(\V U \V_{(W^{1,1}([-L,L]))^*} + \V \T U \V_{L^\infty}), 
\end{align*}
which concludes the proof of Proposition \ref{mainprop}. \end{Proof}

\begin{Corollary}\label{ophasclosedrange} The operators $\T_\A$ and $\T_{\A^C}$ have closed range and finite-dimensional kernel.

\end{Corollary}
\begin{Proof}
Let $\mathcal{R} = \mathcal{I} \circ \mathcal{R}_L$, $X,Y = L^\infty(\R, \C^n)$, $Z = (W^{1,1}([-L,L],\C^n))^*$.  The result then follows for $\T_\A$ from Lemma \ref{AbsClRg} and Proposition \ref{mainprop}.  For $\T_{\A^C}$, let $\mathcal{R}: C^0(\R,\C^n) \to (W^{1,1}([-L,L],\C^n))^*$ be defined analogously, and let $X,Y = C^0(\R, \C^n)$, and the same is true. \end{Proof}

\paragraph{Adjoint properties.} In order to show that the cokernels of $\T_\A$ and $\T_{\A^C}$ are finite-dimensional, we consider the kernels of the adjoint operators $\T_\A^*, \T_{\A^C}^*$. 
Consider first $\T_\A^*: (L^\infty(\R,\C^n))^* \to (L^\infty(\R,\C^n))^*$. 
Abstractly, the adjoint $\T_\A^*$ is defined only as an operator on $(L^\infty(\R,\C^n))^* \cong ((L^\infty(\R,\C))^*)^n$, where $(L^\infty(\R,\C))^*$ can be identified with the space of absolutely continuous finite Borel measures on $\R$. However, we see that for an $n$-tuple $\mu$ of measures in the kernel of $\T_\A^*$, we must have \begin{equation}\label{e:adjker}
    \int U(\xi) \cdot d\mu = - \int \left[\int_\R K(\xi - \xi';\xi)A^{-1}(\xi')U(\xi')d\xi'\right] \cdot d\mu(\xi)
\end{equation} for all $U \in L^\infty$, where $\cdot$ refers here to the dot product on $\C^n$. Note that every component of the matrix-valued function $K(\xi - \xi';\xi)A^{-1}(\xi')$ is an $L^\infty$ function of $\xi$, with $L^\infty$ norm bounded over $\xi'$. Then for any $i,j$, the function $\int (K(\xi - \xi';\xi)A^{-1})_{ij}(\xi') d\mu_i$ is in fact an $L^\infty$ function of $\xi'$. Therefore, since \eqref{e:adjker} must hold for every $U$, we see by equating terms that the $n$-tuple of measures $\mu$ is given by an element of $L^\infty(\R,\C^n)$ through $(\mu)_i = \left(\int (K(\xi - \xi';\xi) A^{-1}(\xi')e_i) \cdot d\mu \right) \lambda$, where $\lambda$ is the Lebesgue measure and $e_i$ is the $i$th standard basis vector in $\C^n$. 

The same is true for $\T_{\A^C}^*$. Since the dual of $C^0(\R,\C)$ can be identified with the space of finite, finitely-additive complex measures on $\R$, we can again identify the dual of $C^0(\R,\C^n)$ with $n$-tuples of measures, and we must have for all $\mu \in (C^0(\R,\C^n))^*$, $U \in C^0$, that \begin{equation*}
    \int U(\xi)\cdot d\mu = - \int \left[\int_\R K(\xi - \xi';\xi)A^{-1}(\xi')U(\xi')d\xi'\right] \cdot d\mu.
\end{equation*}
Now, $K(\cdot - \xi', \cdot) \in C^0(\R, M_n(\C)) \subset W^{1,1}(\R, M_n(\C))$, with $C^0$ norm bounded over $\xi' \in \R$. Therefore, for all $i,j,$ $\int (K(\xi - \xi';\xi)A^{-1}(\xi'))_{ij}d\mu_i$ is well-defined and an $L^\infty$ function of $\xi'$. Then, by the same calculation as above, the $n$-tuple of measures $\mu$ corresponds to an element of $L^\infty(\R, \C^n)$ by $(\mu)_i = \left(\int (K(\xi - \xi';\xi) A^{-1}(\xi')e_i) \cdot d\mu \right) \lambda$. 

For such an $n$-tuple of measures that corresponds to an $L^\infty(\R,\C^n)$ function, the actions of both the operators $\T_\A^*$ and $ \T_{\A^C}^*$ coincide with the action of the operator \[\T^*_{L^\infty}: U(\xi) \mapsto A^*(\xi)U(\xi) + \int K^*(\xi'- \xi; \xi')U(\xi')d\xi'\] on the $L^\infty$ function, where $K^*, A^*$ refer to the conjugate transposes of the matrices $K, A$.

Then for both $\T_\A$ and $\T_{\A^C}$, the kernel of the adjoint operator will be finite-dimensional provided that the kernel of $\T^*_{L^\infty}$ is finite-dimensional on $L^\infty$.

\begin{Lemma}\label{adjoint} The operator $\T^*_{L^\infty}$ satisfies Hypothesis \ref{hyp1} and condition (ii) from Theorem \ref{main} whenever $\T_\A$ (or $\T_{\A^C}$) does. \end{Lemma}  
\begin{Proof} We note that the conditions that $K(\xi -\cdot, \xi) \in W^{1,\infty}(\R, W^{1,1}(\R, M_n(\C))$ and $K(\cdot + \xi, \cdot) \in W^{1,\infty}(\R, W^{1,1}(\R,M_n(\C)))$ exactly guarantee that Hypothesis \ref{hyp1} is readily satisfied by both $\T_\A$ (or $\T_{\A^C}$) and $\T^*_{L^\infty}$. Condition $(ii)(a)$ holds for $\T^*_{L^\infty}$ since $A^*(\xi)$ is invertible whenever $A(\xi)$ is, and we can see that condition $(ii)(b)$ holds by noting that for a hyperbolic $\T_{\A^0}$, we have that $\det((\wh{K^0}(i\ell) + A^0)) = (-1)^n\det((\wh{K^0}(i\ell) + A^0)^*)$, so $\T^\pm$ will be hyperbolic exactly when $(\T^*_{L^\infty})^\pm$ are. \end{Proof}

\paragraph{Limit Operators and the Fredholm Index.} 

The last part of Theorem 1 concerns the indices of the operators, when they are Fredholm. In particular, it will be useful later to have the following fact: 

\begin{Proposition}\label{limitops}
If $\T_\A$ (or $\T_{\A^C}$) is Fredholm, its index depends only on $A^\pm$ and $K^\pm$.
\end{Proposition}

\begin{Proof} Suppose $\T_{\A_1}$ and $\T_{\A_2}$ (or $\T_{\A^C_1}$, $\T_{\A^C_2}$ respectively) satisfy Hypothesis \ref{hyp1} and are Fredholm. We use in the following that, as a consequence of the necessary part of Theorem \ref{main} which is proved below, we may assume that condition (ii) from Theorem \ref{main} is met. Suppose $\A_1, \A_2$ are given by $ (A_1, K_1), (A_2, K_2)$, with the same limits $A^\pm, K^\pm$.

Two Fredholm operators $T_1$ and $T_2$ have the same index if there exists an invertible operator $B$ such that $T_1 - T_2B$ is compact.
The multiplication operator $U(\xi) \mapsto A_2^{-1}(\xi)A_1(\xi)U(\xi)$ is invertible, so we would like to show that $\T_R := \T_{\A_1} - \T_{\A_2}A_2^{-1}A_1$ is compact. We will do this by showing it is the operator-norm limit of a sequence of compact operators. Essentially, we would like to cut off the operator outside $|\xi| < L$, and show that $\T_R$ is the limit of the truncated operators as $L \to \infty$. 

More formally, let $\T_R^L U = E_L \circ \mathcal{I}_L \circ \mathcal{R}_L \circ (\chi^L  \T_R U)$, where
\begin{align*} &\mathcal{R}_L: W^{1,\infty}(\R, \C^n) \to W^{1,\infty}([-L,L],\C^n)) \textrm{ (or }C^1, C^1) \textrm{ is the restriction operator},\\
&\mathcal{I}_L: W^{1,\infty}([-L,L],\C^n) \to L^{\infty}([-L,L],\C^n)\textrm{ (or }C^1, C^0) \textrm{ is the inclusion operator}, \\
&E_L: L^\infty([-L,L],\C^n) \to L^\infty(\R, \C^n)\textrm{ (or }C^0,C^0)\textrm{ extends by 0 outside }[-L,L], \textrm{and }\\
&\chi^L \textrm{ is a smooth characteristic function equal to 1 on }[-(L-1),L-1],\textrm{ and 0 outside } [-L,L].\end{align*}

The operator $\mathcal{I}_L$ is compact, and the operators $\mathcal{R}_L, E_L,$ and multiplication by $\chi^L$ are bounded, so $\T_R^L$ is a compact operator. 

We now show that $\T_R$ is the operator limit of the sequence $\{\T_R^L\}, L \to \infty$. Let $\ep > 0$, and note that $\lv \cdot \rv_{L^\infty}$ may refer equivalently to the supremum norm on $L^\infty$ or $C^0$. We have 
\[
\T_R U(\xi) = (\T_{\A_1} - \T_{\A_2}A_2^{-1}A_1)U(\xi) = \int_\R \left[K_1(\xi - \xi';\xi)U(\xi') - K_2(\xi - \xi';\xi)A_2^{-1}(\xi')A_1(\xi')U(\xi')\right]d\xi'.
\]
Let $A_R(\xi') = \Id - A_2^{-1}(\xi')A_1(\xi')$. Then 
\begin{align*} \lv(\T_R - \T_R^L)U\rv_{L^\infty} &\le \sup_{|\xi| > (L-1)} |\T_R U (\xi)| \\
&\le \sup_{|\xi| > (L-1)}\left( \left| \int_\R (K_1 - K_2)(\xi - \xi';\xi)U(\xi')d\xi'\right| + \left|\int_\R K_2(\xi-\xi';\xi)A_R(\xi')U(\xi')\right| \right) .  \end{align*} 
Since $K_1$ and $K_2$ each converge to $K^\pm$ in $L^1$, we can find $L >0$ that $\sup_{|\xi| > (L-1)} \lv (K_1 - K_2)(\cdot, \xi) \rv_{L^1}  < \frac{\ep}{2}$. Then
\[
\sup_{|\xi| > L-1} \left| \int_\R (K_1 - K_2)(\xi - \xi';\xi)U(\xi')d\xi'\right| \le \sup_{|\xi| > L-1} \lv (K_1 - K_2)(\cdot, \xi) \rv_{L^1} \lv U \rv_{L^\infty} < \frac{\ep}{2} \lv U\rv_{L^\infty}.
\] 
We can then find $M_1$ large enough that $\sup_{\xi' > M_1} \lv A_R(\xi') \rv_{M_n(\C)} < \frac{\ep}{4\sup_\xi\lv K_2(\cdot, \xi) \rv_{L^1}}$, since $A_R(\xi')$ goes to 0 as $|\xi'| \to \infty$. We can also find $M_2 $ large enough that $|\int_{|\zeta| > M_2} K_2(\zeta;\xi)d\zeta| < \frac{\ep}{4\max\left(1, \sup_{\xi'}\lv A_R(\xi')\rv_{M_n(\C)}\right)}$ for all $|\xi| > (L-1)$, since $K_2$ varies continuously in $L^1$ with $\xi$, with limits at infinity. Then take $L> M_1 + M_2 + 1$. 

We have \begin{align*} &\sup_{|\xi| > (L-1)} \left|\int_\R K_2(\xi-\xi';\xi)A_R(\xi')U(\xi')d\xi'\right|\\
     \le& \sup_{|\xi| > (L-1)}\left(\left| \int_{|\xi'| < M_1}K_2(\xi-\xi';\xi)A_R(\xi')U(\xi')d\xi'\right| + \left|\int_{|\xi'| > M_1} K_2(\xi-\xi';\xi)A_R(\xi')U(\xi')d\xi'\right|\right) \\
    \le& \sup_{|\xi| > (L-1)}\left(\sup_{\xi'\in \R}\lv A_R(\xi')\rv \left|\int_{|\zeta| > M_2} K_2(\zeta,\xi)d\zeta\right|  + \sup_{\xi' > M_1} \lv A_R(\xi')\rv \lv  K_2(\cdot, \xi)\rv_{L^1}\right) \lv U \rv_{L^\infty} \\
    < &\left(\frac{\ep}{4} + \frac{\ep}{4}\right)\lv U \rv_{L^\infty}.
\end{align*}

Putting this together, for $L$ sufficiently large, we get \begin{align*}\lv(\T_R - \T_R^L)U\rv_{L^\infty} &\le \sup_{|\xi| > (L-1)}\left( \left| \int_\R (K_1 - K_2)(\xi - \xi';\xi)U(\xi')d\xi'\right| + \left|\int_\R K_2(\xi-\xi';\xi)A_R(\xi')U(\xi')\right| \right) \\ &< \left(\frac{\ep}{2} + \frac{\ep}{2}\right)\lv U \rv_{L^\infty} = \ep\lv U \rv_{L^\infty}.  \end{align*}  So $\T_R$ is the operator-norm limit of the sequence $\{\T_R^L\}$, $L \to \infty$. Then since it is the limit of a sequence of compact operators, $\T_R$ must be compact. Therefore the Fredholm indices of $\T_\A$ and $\T_{\A^C}$ depend only on the limiting operators. 
\end{Proof}

Given this, we are finally ready to prove Theorem 1. 

\begin{Proof}[of Theorem 1, Sufficiency] Assume Hypothesis \ref{hyp1} and conditions (a) and (b) from Theorem \ref{main}, and let $\T$ refer either to $\T_\A$ or $\T_{\A^C}$. From Corollary \ref{ophasclosedrange}, we conclude that $\T$ has closed range and finite-dimensional kernel. 

Using Lemma \ref{adjoint}, $\T^*_{L^\infty}$ then satisfies the same hypotheses as $\T$, meaning the same conclusions apply, and $\T^*_{L^\infty}$ has finite-dimensional kernel. Since the kernel of $\T$ corresponds to the kernel of $\T^*_{L^\infty}$, then $\T$ have finite-dimensional kernel also. 

These together imply that $\T$ is a Fredholm operator.

Lastly, by Proposition \ref{limitops}, the Fredholm index of $\T$ depends only on the limits $A^\pm, K^\pm$. Then the sufficiency part of Theorem 1 is proven. 
\end{Proof}

\begin{Remark}
The above argument extends readily to $L^p$, $1 < p < \infty$, with only minor modifications. In particular, replace $L^\infty$ and $(W^{1,1})^*$ with $L^p$ and $(W^{1,q})^*$ everywhere they appear, with $\frac{1}{q} + \frac{1}{p} = 1 $. Also replace $L^1$ and $W^{1,1}$ with $L^q$ and $W^{1,q}$ in the proof of Proposition \ref{bddH-1}, and the main inequality in the proof of Proposition \ref{bddH-1} with $\lv \wt{\mathcal{K}}U\rv_{W^{1,q}}^q \le C \lv U \rv_{L^q}^q$, where \begin{equation*}C = \sup_{\xi} \lv  K^*( \xi-\cdot, \xi)\rv_{L^1}^{p/q}\sup_\zeta \lv K^*(\cdot - \zeta, \cdot) \rv_{L^1} + \sup_{\xi} \lv  \frac{d}{d\xi}K^*( \xi - \zeta, \xi)\rv_{L^1(\zeta)}^{p/q}\sup_\zeta \lv \frac{d}{d\xi}K^*(\xi - \zeta, \xi) \rv_{L^1(\xi)},
\end{equation*}
which is verified by a short calculation using H\"older's inequality. Lastly, in the proof of Proposition \ref{limitops}, replace $\lv \cdot \rv_{L^\infty}$ with $\lv \cdot \rv_{L^p}^p$ and $\sup_{\xi > (L-1)}(\cdot)$ with $\int_{\xi > (L-1)}|\cdot|^p$ whenever they appear, and note that the inequality $|a+b|^p \le 2^p(|a|^p + |b|^p)$ introduces a factor of $2$. The rest holds without further modification. 

We do note that for $1 < p < \infty$, the requirement that $K^\pm \in L^2_1$ can be dropped, since step 2 in the proof of Proposition \ref{mainprop} can be proved instead using properties of Fourier multipliers. 
\end{Remark}

\subsection{Proof of Theorem 1, Necessity of Conditions for Fredholm}

We now prove the necessity of conditions $(ii)(a)$ and $(ii)(b)$ in Theorem \ref{main}. The proof relies on construction of two Weyl sequences: first, when the principal part has a zero, a sequence becoming concentrated around the zero; second, when the limiting operators are not invertible, a sequence concentrating at the kernel in Fourier space, whose support in physical space is pushed out to infinity. These sequences are constructed for $L^\infty$ and $C^0$, but can be easily modified for the $L^p$ case. 


\paragraph{Construction of Weyl Sequence for condition (ii)(b).} Assume that condition $(ii)(b)$ is not met; that is, $A^{-1} \notin L^\infty(\R, M_n(\C))$. We will construct a Weyl sequence first for $\T_\A$, then for $\T_{\A^C}$. 

First, consider $\T_\A$, so that $A(\cdot)$ is in $L^\infty$ but not necessarily $C^0$. Since $A^{-1} \notin L^\infty$, then for each $N \in \N$ there exists a set $E_N \subset \R$ of positive measure such that $\inf_{|v| = 1} | A(\xi) v | < \frac{1}{2N}$. If the measure of $E_N$ is greater than $\frac{1}{N}$, redefine it as a subset so that its measure is equal to $\frac{1}{N}$. Since $A(\xi) \in L^\infty(\R, M_n(\C)),$ by Lusin's theorem, there exists for each $N \in \N$ an $A_N \in C^0(\R, M_n(\C)),$ and $ \Omega_N \subset \R$, so that $A_N(\xi) \equiv A(\xi)$ on $\Omega_N$, with the measure of $\R \ \backslash \ \Omega_N$ equal to half the measure of $E_N$. Let $\wt{E}_N = E_N \cap \Omega_N$. Because each $A_N$ is continuous, we can define $v_N(\xi)$ supported on $\wt{E}_N$, piecewise constant, so that $| A(\xi)v_N(\xi)| < \frac{1}{N}$, as follows: 

$A_N$ is uniformly continuous on each $[k,k+1] \subset \R, k \in \N$. Cover $[k,k+1]\cap \wt{E}_n$ with disjoint intervals $\Delta_{k,i}$ of length at most $\delta_k$ such that $| (A_N(\xi) - A_N(\xi_0))v | < \frac{1}{2N}$ for any $|\xi - \xi_0| < \delta_k, |v| = 1$. Choose $\xi_{k,i} \in \Delta_{k,i} \cap \wt{E}_N$, and let $v_{k,i}, |v_{k,i}| = 1$ be such that $|A_N(\xi_{k,i})v_{k,i}| < \frac{1}{2N}$. Let
\[
u_N(\xi) = \begin{cases} v_{k,i}, \qquad &\xi \in \Delta_{k,i} \cap \wt{E}_N \\ 0, \qquad &\textrm{otherwise}. \end{cases}
\]
Then $u_N$ is a measurable function, with $\lv u_N \rv_{L^\infty} = 1$, for which $|A(\xi)v_N(\xi)| < \frac{1}{N}$ for all $\xi$. 

For $\T_{\A^C}$ the construction is somewhat simpler. Since in this case $A(\cdot)$ is continuous with $A^{-1} \notin L^\infty$, there exists for all $N \in \N$ an interval $E_N$ with positive measure so that $\inf_{|v| = 1}|A(\xi)v| < \frac{1}{2N}$ for all $\xi \in E_n$. Again, redefine $E_N$ possibly as a subinterval so that $m(E_N) < \frac{1}{N}$. Let $\xi_N \in E_N$, and let $v_0, |v_0| = 1$, be a vector such that $|A(\xi_N)v_0| < \frac{1}{2N}$. Then there exists a subinterval $\wt{E}_N$ of $E_N$ with positive measure such that $|A(\xi)v_0| < \frac{1}{N}$ for $\xi \in \wt{E}_N$. Let $\chi_N$ be a smooth function supported on $\wt{E}_N$ with $\lv \chi_N \rv_{C^0} = 1.$ Then let
\[
u_N(\xi) = \chi_N(\xi)v_0.\] 
Note that in both constructions, we get that $u_N$ is supported on a set of measure $\wt{\delta}_N \le \frac{1}{N}$. 
%
%
%
\begin{Proof}[(Necessity of Condition (ii)(b))] We will now prove that $\{u_N\}$ is a Weyl sequence. Let $\T$ refer either to $\T_\A$ or $\T_{\A^C}$, and let $\{u_N\}$ be the corresponding sequence defined above. Let $\lv \cdot \rv_{L^\infty}$ refer equivalently to the norm on $L^\infty$ or $C^0$. Since $\lv A(\xi) u_N(\xi) \rv_{L^\infty} < \frac{1}{N}$, we can choose $N_0$ large enough so that $\lv A(\xi) u_N (\xi) \rv_{L^\infty} < \frac{\ep}{2}$ for all $N > N_0$. Note also that $\lv u_N \rv_{L^1} \le \wt{\delta_N} \to 0$, and, by the assumptions on $K$, that $\sup_{\xi,\xi'}\lv K(\xi - \xi',\xi)\rv_{M_n(\C)} < \infty$. Then we can find $N_0$, possibly larger, so that for $N > N_0,$ 
\[
\lv \int_\R K(\xi - \xi', \xi)u_N(\xi')d\xi' \rv_{L^\infty} \le \sup_{\xi,\xi'}\lv K(\xi - \xi',\xi)\rv_{M_n(\C)} \lv u_N \rv_{L^1}< \frac{\ep}{2}
.\] 
Putting this together, there exists $N_0$ large enough so that for $N > N_0$,

\begin{align*} \lv \T u_N \rv_{L^\infty} &\le \lv A u_N \rv_{L^\infty} + \lv \int_\R K(\xi - \xi', \xi)u_N(\xi)d\xi' \rv_{L^\infty} \\
&< \frac{\ep}{2} + \frac{\ep}{2} = \ep.
\end{align*}

Then $\{u_N\}$ forms a Weyl sequence for $\T$, and $\T$ is not Fredholm, showing the necessity of condition (ii)(a).  \end{Proof}
%
%
%

\paragraph{Construction of Weyl sequence for condition (ii)(a).} Next, assume that Condition (ii)(a) is not satisfied; i.e., there exists $m \in \R$ such that \[ \det(\wh{K^\pm}(i m) + A^\pm) = 0. \] We construct one Weyl sequence for both $\T_\A$ and $\T_{\A^C}$. Without loss of generality suppose $\det(\wh{K^+}(i m) + A^+) = 0$. Then there exists a vector $v \in \R^n, | v | = 1$, so that $(\wh{K^+}(i m) + A^+)v = 0$. Let
%
%
\begin{align*}u(\xi) &= e^{-\frac{1}{4}\xi^2 - im\xi}v,\\ u_N(\xi) &= u(\frac{\xi - N^2}{N}).\end{align*}
Note that $\ds \widehat{u_N}(i \ell) = N\sqrt{2}e^{-N^2(\ell-m)^2}e^{-iN^2}v$, and $\lv \wh{u}_N \rv_{L^1} = \sqrt{2\pi}$ for all $N$.

\begin{Proof}[(Necessity of condition (ii)(a))] We now prove $\{u_N\}$ is a Weyl sequence. Again, let
$\T$ refer to either $\T_\A$ or $\T_{\A^C}$, and $\lv \cdot \rv_{L^\infty}$ to the norm on $L^\infty$ or $C^0$ equivalently. We have \begin{align*} 
 \left\V \T u_N \right\V_{L^\infty} &= \lv A(\xi)u_N(\xi) + \int_{\R}K(\xi - \xi'; \xi)u_N(\xi')d\xi')\rv_{L^\infty}  \\
 &\le \left\V (A(\xi) - A^+)u_N\right\V_{L^\infty} + \left\V  \int_{\R}(K(\xi - \xi'; \xi) - K^+(\xi - \xi'))u_N(\xi')d\xi' \right\V_{L^\infty}  + \left\V \T^+ u_N \right\V_{L^\infty}. 
\end{align*}
Let $\ep > 0$. Choose $N_0$ large enough that for $N > N_0$, \[\ds \sup_{\xi \le n} |u_N(\xi)| < \frac{\ep}{\max(1,\sup_{\xi \le n}\lv A(\xi) - A^+\rv^2)},\] which we may do by the choice of $u_N$, and so that $ \sup_{\xi > N}\lv A(\xi) - A^+\rv <  \ep$, which we may do by Hypothesis \ref{hyp1}. Then for $N > N_0$, \begin{align*}  \lv(A(\xi) - A^+)u_N(\xi)\rv_{L^\infty} &= \max\left(\sup_{\xi \le n} \left|(A(\xi) - A^+)u_N(\xi)\right|,  \sup_{\xi > n} \left|(A(\xi) - A^+)u_N(\xi)\right|\right) \\
& \le \max\left(\sup_{\xi \le n}\lv A(\xi) - A^+\rv \sup_{\xi \le n}|u_N|,  \sup_{\xi > n}\lv A(\xi) - A^+\rv \sup_{\xi > n}|u_N|\right)\\
& < \ep.
\end{align*}



Next, choose $M$ large enough so that \[\ds \sup_{\xi >M} \lv K(\cdot, \xi) - K^+(\cdot) \rv_{L^1} < \ep.\]

Note that  $\ds \lv\int_{\R} \left[K(\xi - \xi'; \xi) - K^+(\xi - \xi')\right]u_N(\xi')d\xi' \rv_{L^\infty} $ is equal to the larger of 

\[\ds \sup_{\xi \le M} \left|\int_{\R}  \left[K(\xi - \xi'; \xi) - K^+(\xi - \xi')\right]u_N(\xi')d\xi' \right| , \sup_{\xi > M} \left|\int_{\R} \left[K(\xi - \xi'; \xi) - K^+(\xi - \xi')\right]u_N(\xi')d\xi' \right|;\]

and by the choice of $M$, 

\[\ds \sup_{\xi > M} \left |\int_{\R} \left[K(\xi - \xi'; \xi) - K^+(\xi - \xi')\right]u_N(\xi')d\xi'  \right| \le\sup_{\xi>M} \lv  K(\cdot, \xi) - K^+(\cdot) \rv_{L^1}\lv u_N \rv_{L^\infty} < \frac{\ep}{4}.\]

On the other hand, we have

\begin{align*} \ds &\sup_{\xi \le M} \left|\int_{\R}  \left[K(\xi - \xi'; \xi) - K^+(\xi - \xi')\right]u_N(\xi')d\xi' \right|\\
 = &\sup_{\xi \le M} \left|\int_{\xi' \le \xi + L}  \left[K(\xi -\xi',\xi) - K^+( \xi - \xi')\right]{u_N}(\xi')d\xi'+\int_{\xi' >\xi + L}  \left[K(\xi -\xi',\xi) - K^+( \xi - \xi')\right]{u_N}(\xi')d\xi' \right|. \end{align*}

Since $K(\cdot, \xi) \in W^{1,\infty}(\R, W^{1,1}(M_n(\C)))$, with limits at infinity, we can find $L(\ep)$ so that 
$$
\int_{|\zeta| > L}|K(\zeta,\xi) - K^+(\zeta)|d\zeta' < \frac{\ep}{2} \textrm{, for all } \xi,
$$
yielding
\[
\sup_{\xi \le M} |\int_{\xi' >\xi + L}  \left[K(\xi -\xi',\xi) - K^+( \xi - \xi')\right]{u_N}(\xi')d\xi' | \le \left(\sup_{\xi} \int_{-\infty}^{-L} | K(\xi', \xi) - K^+(\xi')| d\xi'\right)\lv u_N \rv_{L^\infty} \le \frac{\ep}{2}
.\]
Also, because $u_N$ is shifted to the right by $N^2$, but is only stretched by a factor of $N$, we can choose $N_0$ large enough so that for $N > N_0$, $\ds \sup_{\xi \le L + M} \left| u_N(\xi) \right| < \frac{\ep}{2\sup_{\xi} \lv  K(\cdot, \xi) - K^+(\cdot) \rv_{L^1}}$, yielding
\[
\sup_{\xi \le M} \int_{\xi' \le \xi L}  \left[K(\xi -\xi',\xi) - K^+( \xi - \xi')\right]{u_N}(\xi')d\xi' \le \sup_{\xi} \lv  K(\cdot, \xi) - K^+(\cdot) \rv_{L^1}\sup_{\xi < L+M} |u_N(\xi)| \le \frac{\ep}{2}.
\]
This, combined with the above, gives
\[
\sup_{\xi \le M} \left|\int_{\R}  \left[K(\xi - \xi'; \xi) - K^+(\xi - \xi')\right]u_N(\xi')d\xi' \right| <\frac{\ep}{2}+\frac{\ep}{2} = \ep,
\]

so that \[\ds \lv\int_{\R} \left[K(\xi - \xi'; \xi) - K^+(\xi - \xi')\right]u_N(\xi')d\xi' \rv_{L^\infty} < \ep.\]

Lastly, for any $\ep > 0$ we can choose $N_0$ large enough so that $  \lv A^+ + \wh{K^+}(i\ell)\rv_{M_n(\C)} < \frac{\ep}{2}$ for $|m - \ell| < \frac{1}{\sqrt{N_0}}$, 

and also large enough so that for $N > N_0$, $\ds \int_{|m - \ell| \ge \frac{1}{\sqrt{n}}} \left| \wh{u_N}(i\ell)\right| d\ell < \frac{\ep}{2 \sup_{\nu \in \R}\lv(A^+ + \wh{K^+}(i\nu))\rv}$, because the $\wt{u}_N$ become increasingly localized about $\ell = m$. This gives that  \begin{align*} \lv \T^+ u_N \rv_{L^\infty} &\le \frac{1}{\sqrt{2\pi}}\lv \wh{\T^+ u_N} \rv_{L^1} \\ &= \frac{1}{\sqrt{2\pi}}\lv \left( A^+ + \wh{K^+}(i\ell)\right) \wh{u_N} \rv_{L^1} \\
 &= \int_{|m - \ell| \ge \frac{1}{\sqrt{n}}} \left| \left( A^+ + \wh{K^+}(i\ell)\right) \frac{\wh{u_N}(i\ell)}{\sqrt{2\pi}}\right|d\ell + \int_{|m - \ell| < \frac{1}{\sqrt{n}}} \left|\left( A^+ + \wh{K^+}(i\ell)\right) \frac{\wh{u_N}(i\ell)}{\sqrt{2\pi}}\right| d\ell \\
& \le \sup_{\nu \in \R}\lv \left(A^+ + \wh{K^+}(i\nu)\right)\rv_{M_n(\C)} \int_{|m - \ell| \ge \frac{1}{\sqrt{n}}} \left|\wh{u_N}(i\ell)\right| d\ell +   \frac{\ep}{2} \int_{|m - \ell| < \frac{1}{\sqrt{n}}} \left| \frac{\wh{u_N}(i\ell)}{\sqrt{2\pi}}\right|d\ell \\
& < \frac{\ep}{2} + \frac{\ep}{2} = \ep. \end{align*}

Therefore, for any $\ep >0$, for $N > N_0$ large enough that the previous inequalities hold, we get that

\begin{align*} \left\V \T u_N \right\V_{L^\infty} &\le \left\V (A(\xi) - A^+)u_N\right\V_{L^\infty} + \left\V  \int_{\R}(K(\xi - \xi'; \xi) - K^+(\xi - \xi'))u_N(\xi')d\xi' \right\V_{L^\infty}  + \left\V \T^+ u_N \right\V_{L^\infty}\\
&\le \ep + \ep + \ep = 3\ep,
 \end{align*}
 which implies that $$ \lim_{N \to \infty} \left\V \T u_N \right\V_{{L^\infty}} = 0, \textrm{ with } \lv u_N \rv_{L^\infty} = 1.$$

Thus $\{u_N\}$ is a Weyl sequence for $\T$, which implies that $\T$ is not Fredholm, showing the necessity of Condition (ii)(b). \end{Proof}


\section{Spectral Flow and the Fredholm Index}\label{s:3}

We establish results that allow us to compute the index of the nonlocal operator 
\[ 
\T: U(\xi) \mapsto A(\xi)U(\xi) + \int_{\R}K(\xi - \xi'; \xi)U(\xi')d\xi',
\]
defined in Section \ref{s:2.1}, in many specific situations. Assuming exponential localization of convolution kernels, stronger than in Section \ref{s:2}, the Fredholm index is given by the \it spectral flow\rm  \ of an operator with the same limits at infinity. The approach here is somewhat closely following  \cite{faye2014fredholm}, which in turn is relying on ideas from \cite{robbin1995spectral,mallet1999fredholm}. The argument in \cite{faye2014fredholm} needs to be modified for two reasons: a change in the form of the operator, and a change in the domain of the operator (from $L^2$ to $L^\infty, C^0$). The latter affects the argument only in Lemma \ref{forsimple}, and the former is the cause of the rest of the modifications. Technically, one needs to carefully inspect the space of allowed perturbations so that spectral crossings along relevant paths are generic.  

In order to state and prove the following theorem, we require two additional assumptions on the convolution kernel which were not needed in Sections \ref{s:2.1}-\ref{s:2.2}.

\begin{Hypothesis}\label{hyp3.1} The generalized convolution kernel $K$ is exponentially localized in its first argument---that is, for some $\eta >0$, we have $K \in C^0(\R, W_\eta^{1,1}(\R, M_n(\C))),$ where \[W_\eta^{1,1}(\R, M_n(\C)) = \left\{ f \in W^{1,1}(\R, M_n(\C)) \ \ \Big| \ \max_{1\leq j,k\leq n}\left(\lv f_{j,k}(\cdot)e^{\eta |\cdot|}\rv_{L^1} + \lv f_{j,k}'(\cdot)e^{\eta |\cdot|}\rv_{L^1}\right)< \infty\right\}.\] \end{Hypothesis}

\begin{Hypothesis}\label{hyp3.2} The Fourier transforms \[\nu \mapsto \wh{K}^\pm(\nu) + A^\pm\] extend to bounded analytic functions in the strip $S_\eta = \{ \nu \in \C \ | \ |\Re(\nu)|< \eta\}$. \end{Hypothesis}

\begin{Theorem}\label{indexflow} Let $\T$ refer either to $\T_\A$ or $\T_{\A^C}$ as defined in \eqref{e:opdef}. Suppose Hypotheses \ref{hyp1}, \ref{hyp3.1}, and \ref{hyp3.2} are satisfied, as well as condition (ii) of Theorem \ref{main}. Let $A^C(\cdot)$ be a continuous function, possibly different from $A(\cdot)$ if the latter is not continuous, with $(A^C)^{-1}(\cdot) \in L^\infty(\R, M_n(\C))$, such that $\lim_{\xi \to \pm\infty}A^C(\xi) = A^\pm$. Suppose that for the operator $\T^C$ defined by $A^C$ and $K$, there exist only finitely many values $\xi_0 \in \R$ for which $\T^C$ is not hyperbolic; that is, for which $\det(\wh{K}_{\xi_0}(i\ell) + A^C(\xi_0)) = 0$ for some $\ell \in \R$. 

Then the Fredholm index of $\T$ is given by 
\[ 
\mathrm{ind}\,\T = -\textrm{cross}(\A),
\]
where $\textrm{cross}(\A)$ denotes  the net number of roots, counted with multiplicity, of the characteristic equation 
\begin{equation}\label{char} 
d^\xi(\nu) = \det\left(\wh{K}_\xi(\nu)+ A^C(\xi)\right)
\end{equation}
that cross the imaginary axis from left to right as $\xi$ is increased from $-\infty$ to $+\infty$; see \eqref{crossnum} below for a more precise definition. 
\end{Theorem}

The remainder of this section will be devoted to more precisely stating and proving Theorem \ref{indexflow}. 
%
In particular, we will prove the following theorem, from which Theorem \ref{indexflow} follows. For notational simplicity, in the following we identify the symbol $\A$ with its associated operators $\T_\A$ and $\T_{\A^C}$, suppressing the difference in domains. Because the indices of $\T_\A$, $\T_{\A^C}$ depend only on the limits $A^\pm, K^\pm(\cdot)$, we denote the Fredholm index of $\T_\A$ by $\iota(\A^-, \A^+)$, and the Fredholm index of $\T_{\A^C}$ by $\iota_c(\A^+, \A^-)$. 
We also define, for a constant-coefficient operator $\A^0= (A^0, K^0(\cdot))$ the function \begin{equation}\label{e:delta} \Delta_{\A^0}(\nu) = \wh{K}^0(\nu) + A^0 \end{equation} and the characteristic equation \begin{equation} d^0(\nu) = \det(\Delta_{\A^0}) = 0. \end{equation}

\begin{Theorem}\label{arho} Let $\{\A^\rho\}$, for $\rho \in \R$, be a continuously varying one-parameter family of constant-coefficient operators $(A^\rho, K^\rho)$, with limit operators $\A^\pm = \lim_{\rho \to \pm \infty}\A^\rho$. We suppose that:
\begin{enumerate}[(i)]
    \item the limit operators $\A^\pm$ are hyperbolic in the sense that for all $\ell \in \R$, \[d^\pm(i\ell) = \det\left(\wh{K}^\pm(i\ell) + A^\pm\right) \neq 0,\]
    
    \item $\Delta_{\A^\rho}$ defined in \eqref{e:delta} is a bounded analytic function in the strip $S_\eta = \{ \lambda \in \C \ \big|\  |\Re(\lambda)| < \eta \}$ for each $\rho \in \R$, and 
    
    \item there are finitely many values of $\rho$ for which $\A^\rho$ is not hyperbolic.
\end{enumerate}
Then 
\[ 
\iota(\A^-,\A^+) = \iota_c(\A^-, \A^+) = -\textrm{cross}(\{\A^\rho\})
\]
is the net number of roots of $d^\rho(\nu) = 0$, counted with multiplicity, which cross the imaginary axis from left to right as $\rho$ is increased from $-\infty$ to $+\infty$; again, see \eqref{crossnum} below for a more precise definition. 
\end{Theorem}
In the proof, we approximate the family $\{\A^\rho\}$ by a generic family. 
To do so, we introduce the set $\mathcal{P} := \mathcal{P}(\R, W^{1,1}_\eta(\R, M_n(\C))\times M_n(\C))$, the Banach space of continuous paths for which conditions (i) and (ii) of Theorem \ref{arho} are satisfied. We also consider the dense set $\mathcal{P}^1 := C^1\left(\R, W^{1,1}_\eta(\R,M_n(\C))\times M_n(\C)\right)\cap \mathcal{P}$. We then first prove that the set of paths with only simple crossings is dense in $\mathcal{P}$. Then, using the proof that for a map with only simple crossings, the Fredholm index is given by the crossing number, the result will follow. 

\paragraph{Notation and Definitions.} 
For a continuously varying one-parameter family $\{\A^\rho\}$ of constant-coefficient operators, a \bf crossing \rm for $\{\A^\rho\}$ is a real number $\rho_j$ for which $\A^{\rho_j}$ is not hyperbolic. The set 
\[ 
NH(\{\A^\rho\}) = \{\rho \in \R \ | \  \textrm{the constant-coefficient operator } {\A^\rho} \textrm{ is not hyperbolic} \}
\]
is the set of all crossings for $\{\A^\rho\}$. Condition (iii) in Theorem \ref{arho} is satisfied only if $NH(\{\A^\rho\})$ is a finite set, which we then can write as $\{\rho_1, ... ,\rho_m\}$. We also have that for any $\{\A^\rho\}$ satisfying the conditions of Theorem \ref{arho} and for any $\rho_j \in NH(\{\A^\rho\})$, the equation \begin{equation} 
d_{\rho_j} := \det(\Delta_{\A^{\rho_j}}) = 0 
\end{equation} 
has finitely many roots in the strip $S_\eta$, due to the analyticity and boundedness of $\Delta_{A^\rho}$, and due to the fact that $d_{\rho_j}(i\ell) \xrightarrow{|\ell| \to \infty} \det(A^{\rho_j}) \neq 0$. Then the \bf crossing number\rm, cross$(\{\A^\rho\})$, can be defined as the net number of roots which cross the imaginary axis as $\rho$ goes from $-\infty$ to $+\infty$, as follows. 

Fix any $\rho_j \in NH(\{\A^\rho\})$ and let $\{\nu_{j,l}\}_{l=1}^{k_j}$ denote the roots of $d_{\rho_j}(\nu)$ on the imaginary axis, listing multiple roots repeatedly according to their multiplicity. Let $M_j$ be the sum of their multiplicities. For $\rho$ near $\rho_j$, with $\pm(\rho-\rho_j) > 0$, this equation has exactly $M_j$ roots, counting multiplicity, near the imaginary axis, $M_j^{L_\pm}$ with negative real part and $M_j^{R_\pm}$ with positive real part. Then the crossing number is defined as 

\begin{equation}\label{crossnum} \textrm{cross}(\A) = \sum_{j=1}^m \left(M_j^{R_+}-M_j^{R_-}\right). \end{equation}

For $\{\A^\rho\} \in \mathcal{P}^1$, a crossing $\rho_j$ is \bf simple \rm if there is exactly one simple root $\nu_*$ of $d_{\rho_j}$ on the imaginary axis, which crosses the imaginary axis with nonvanishing speed as $\rho$ goes through $\rho_j$. For such a crossing, the root can be locally continued as a function of $\rho$, giving a function $\nu(\rho) \in C^1(\R, \C).$ Non-vanishing speed of crossing then corresponds to $\Re(\overset{.}{\nu}(\rho_j)) \neq 0$. 

For a path in $\mathcal{P}^1$ with only simple crossings,  let $\nu_j(\rho)$ be the function defined near  a crossing $\rho_j$ such that $d_\rho(\nu_j) = 0$ and $\Re(\nu_j(\rho_j)) = 0$. Then we have
\[
\textrm{cross}(\{\A^\rho\}) = \sum_{j=1}^m\Re(\overset{.}{\nu_j}(\rho_j)). 
\]
We next prove that the set of paths with simple crossings is dense in $\mathcal{P}$. 

\begin{Lemma}\label{simpledense} Let $\{\A^\rho\} \in \mathcal{P}$, with limit operators $\A^\pm = \lim_{\rho \to \pm\infty}\A^\rho$, such that $NH(\A)$ is a finite set. Then, for $\ep > 0$, there exists $\{\wt{\A}^\rho\} \in \mathcal{P}^1$ such that 
\begin{enumerate}[(i)]
\item $\wt{\A}^\pm = \A^\pm$,
\item $\lv \wt{\A^\rho} - \A^\rho\rv_{W^{1,1}_\eta \times M_n(\C)} < \ep$ for all $\rho \in \R$, and
\item $\{\wt{\A}^\rho\}$ has only simple crossings.
\end{enumerate}
\end{Lemma}
\begin{Remark} If $\ep$ is chosen small enough in the above lemma, then cross($\{\A^\rho\}$) = cross($\{\wt{\A}^\rho\})$,  since the roots of $d^\rho$, which is a holomorphic function, vary continuously in the Hausdorff topology. 
\end{Remark}

In order to prove Lemma \ref{simpledense}, we define submanifolds of $M_n(\C)$. For $0 \le k \le n$, the sets $\mathbf{G}_k \subset M_n(\C)$ and $\mathbf{H} \subset M_n(\C) \times M_n(\C)$ are given by
\begin{align*}
    \mathbf{G}_k& = \{M \in M_n(\C) \ \big| \ \textrm{rank}(M) = k\},\\
    \mathbf{H} &= \{(M_1,M_2) \in (M_n(\C))^2 \ \big| \ \textrm{rank}(M_1) = n-1, M_2 \textrm{ is invertible, and rank}(M_1M_2^{-1}M_1) = n-2\}.
\end{align*}

The sets $\mathbf{G}_k$ and $\mathbf{H}$ are analytic submanifolds of $M_n(\C)$ and $(M_n(\C))^2$ respectively, of complex dimension 
\[
\dim_\C(\mathbb{G}_k) = n^2 - (n-k)^2, \dim_\C(\mathbb{H}) = 2n^2-2; 
\]
see \cite{mallet1999fredholm}. 

For an operator $\A = (A, K(\cdot))$, we rewrite its convolution kernel $K$ more generally as $K(\xi) + B_1\delta_s(\xi - \xi_1) + B_2 \delta_s(\xi - \xi_2)$, where $B_1,B_2$ are real matrices, $\delta_s(\cdot) = \frac{1}{\sqrt{\pi}}e^{-|\cdot|^2}$, and $\xi_1, \xi_2$ are fixed positive real numbers such that $\xi_1/\xi_2$ is irrational.  For an operator of the form considered in Lemma \ref{simpledense}, $B_1, B_2 = 0$.

We then consider the following maps:
\begin{align*}
\mathcal{F},\mathcal{G} &: (W^{1,1}_\eta(\R, M_n(\C)) \times (M_n(\C))^3 ) \times \R \to M_n(\C) \\
\mathcal{F} \times \mathcal{G} &: (W^{1,1}_\eta(\R, M_n(\C)) \times (M_n(\C))^3 ) \times \R \to M_n(\C) \times M_n(\C)\\
\mathcal{D} &: (W^{1,1}_\eta(\R, M_n(\C)) \times (M_n(\C))^3 ) \times T \to M_n(\C) \times M_n(\C)
\end{align*}
given by 
\begin{align*}
\mathcal{F}(\A,\ell) &= \wh{K}(i\ell) + A + B_1\wh{\delta}_s e^{-i\ell\xi_1} + B_2\wh{\delta}_s e^{-i\ell\xi_2} \\
\mathcal{G}(\A,\ell) &= \wh{K}'(i\ell) - B_1e^{-i\ell\xi_1}\left(\xi_1\wh{\delta}_s(i\ell) + \wh{\delta}'_s(i\ell)\right) - B_2e^{-i\ell\xi_2}\left(\xi_2\wh{\delta}_s(i\ell) + \wh{\delta}'_s(i\ell)\right) \\
\mathcal{F} \times \mathcal{G}(\A, \ell) &= (\mathcal{F}(\A,\ell),\mathcal{G}(\A,\ell))\\
\mathcal{D}(\A,\ell) &=(\mathcal{F}(\A,\ell_1),\mathcal{F}(\A,\ell_2))
,\end{align*}
where $T$ is the set \[T = \{(\ell_1,\ell_2) \in \R^2 \ | \ \ell_1 < \ell_2\}.\]

\begin{Proposition} Suppose that $\A = (A,K) \in M_n(\C) \times W^{1,1}_\eta(\R, M_n(\C))$ satisfies the conditions 
\begin{fleqn}\label{conditions}
\begin{align*}
    \hspace{.3 in} &(i)\ \mathcal{F}(\A,\ell) \notin \mathbf{G}_k, &&0 \le k \le n -2, \ell \in \R,\\
   \hspace{.3 in} &(ii)\ (\mathcal{F}\times\mathcal{G})(\A,\ell) \notin \mathbf{G}_{n-1} \times \mathbf{G}_k, &&0 \le k \le n-1, \ell \in \R,\\
   \hspace{.3 in} &(iii) \ (\mathcal{F}\times\mathcal{G})(\A,\ell) \notin \mathbf{H}, &&\ell \in \R, \\
    \hspace{.3 in}&(iv) \ \mathcal{D}(\A, \ell_1,\ell_2) \notin \mathbf{G}_{n-1} \times \mathbf{G}_{n-1},  &&(\ell_1,\ell_2) \in T,
\end{align*}
\end{fleqn}
for all ranges of $k,\ell, \ell_1$, and $\ell_2$. Then the constant-coefficient operator \eqref{coco} has at most one $\ell \in \R$ such that $i\ell$ is a root of $\det \Delta_{\A}(\nu) = 0$, and the root is simple. 

\end{Proposition}

\begin{Proof}
We omit the proof here, as it is identical to \cite[Prop. 4.3]{faye2014fredholm}.
\end{Proof}

\begin{Proposition}\label{surjective} The maps $\mathcal{F}, \mathcal{F}\times \mathcal{G}$, and $\mathcal{D}$ have surjective derivative with respect to the first argument $\A$ at each point $(\A, \ell) \in W^{1,1}_\eta(\R,M_n(\C)) \times (M_n(\C))^3 \times \R $ and $W^{1,1}_\eta(\R,M_n(\C)) \times (M_n(\C))^3 \times T $, respectively. 
\end{Proposition}

\begin{Proof}
From its definition, we see that the derivative of $\mathcal{F}$ with respect to $A$ is $\mathbb{I}_n$, which is surjective onto $M_n(\C)$, and the derivative of $\mathcal{F}\times\mathcal{G}$ with respect to $(A, B_1, B_2)$ is given by the matrix 
\begin{equation}\label{e:mat}
\begin{pmatrix}\mathbb{I}_n &\wh{\delta}_s(i\ell) e^{-i\ell\xi_1}\mathbb{I}_n & \wh{\delta}_s(i\ell) e^{-i\ell\xi_2}\mathbb{I}_n\\ 0_n & e^{-i\ell\xi_1}(\xi_1\wh{\delta}_s(i\ell) + \wh{\delta}'_s(i\ell))\mathbb{I}_n& e^{-i\ell\xi_2}(\xi_2\wh{\delta}_s(i\ell) + \wh{\delta}'_s(i\ell))\mathbb{I}_n\end{pmatrix}. 
\end{equation}
Because $(\xi_1\wh{\delta}_s(i\ell) + \wh{\delta}'_s(i\ell))$ and $(\xi_2\wh{\delta}_s(i\ell) + \wh{\delta}'_s(i\ell))$ are never both equal to 0 at the same value of $\ell$, this operator in \eqref{e:mat} is surjective,  onto $M_n(\C) \times M_n(\C)$.

Now, fixing $(\ell_1, \ell_2) \in T$, we will have that one of the quantities $\xi_1(\ell_1 - \ell_2)$ or $\xi_2(\ell_1 - \ell_2)$ is not a multiple of $2\pi$. Supposing without loss of generality that it is $\xi_1(\ell_1 - \ell_2)$, then we have that the derivative of $\mathcal{D}$ with respect to $(A, B_1)$ is 
\begin{equation*} 
\begin{pmatrix} \mathbb{I}_n & \wh{\delta}_s(i\ell_1)e^{-i\ell_1 \xi_1}\mathbb{I}_n \\ \mathbb{I}_n & \wh{\delta}_s(i\ell_2)e^{-i\ell_2\xi_1}\mathbb{I}_n \end{pmatrix}, \end{equation*}
which is also surjective onto $M_n(\C) \times M_n(\C)$. \end{Proof}

In order to complete the proof of Lemma \ref{simpledense}, we use the notion of transversality for smooth manifolds. A smooth map $f: \mathcal{X} \to \mathcal{Y}$ from two manifolds is transverse to a submanifold $\mathcal{Z} \subset \mathcal{Y}$ on a subset $\mathcal{S}\subset\mathcal{X}$ if 
\[
\textrm{rg}(Df(x)) + T_{f(x)}\mathcal{Z} = T_{f(x)}\mathcal{Y} \textrm{ whenever } x \in \mathcal{S} \textrm{ and } f(x) \in \mathcal{Z},
\]
where $T_p(M)$ denotes the tangent space of $M$ at a point $p$. 

\begin{Theorem}[Transversality Density Theorem]\label{transversality} 
Let $\mathcal{V,X,Y}$ be $C^r$ manifolds, $\Psi: \mathcal{V} \to C^r(\mathcal{X,Y})$ a representation, and $\mathcal{Z} \subset \mathcal{Y}$ a submanifold and $ev_\Psi: \mathcal{V} \times \mathcal{X} \to \mathcal{Y}$ the evaluation map. Assume that: 
\begin{enumerate}[(i)]
\item $\mathcal{X}$ has finite dimension $N$ and $\mathcal{Z}$ has finite codimension $Q$ in $\mathcal{Y}$;
\item $\mathcal{V}$ and $\mathcal{X} $ are second countable;
\item $r > \max(0, N-Q)$;
\item $ev_\Psi$ is transverse to $\mathcal{Z}$. 
\end{enumerate}
Then the set $\{V \in \mathcal{V} \ | \ \Psi_V \textrm{ is transverse to } \mathcal{Z} \}$ is residual in $\mathcal{V}$. 
\end{Theorem}

The proof of this theorem can be found in \cite{abraham1967transversal}. 

\begin{Proposition}\label{resid} There exists a residual (and hence dense) subset of $\mathcal{P}^1$ such that for any $\{\mathcal{A}^\rho\}$ in this subset, all conditions from Proposition \ref{conditions} are satisfied for each $\A^\rho$, $\rho \in \R$. 

\end{Proposition}

\begin{Proof}
We apply the Transversality Density Theorem \ref{transversality} to show that there is a residual subset of $\mathcal{P}^1$ such that all the maps $\mathcal{F, (F \times G), D}$ are transverse to the manifolds appearing in Proposition \ref{conditions} on $(\rho, \ell) \in \R^2$ and $(\rho, \ell_1, \ell_2) \in \R \times T$, respectively. We show the proof for $\mathcal{F}$, the others being similar.

We let $\mathcal{V} = \mathcal{P}^1$, $\mathcal{X} = \R^2$ and $\mathcal{Y} = M_n(\C)$, with submanifold $\mathcal{Z} = \mathbf{G}_k$, for $0 \le k \le n-2$, in the hypotheses of Theorem \ref{transversality}. Then for $\{\A^\rho\} \in \mathcal{P}^1$, we let $\Psi_{\{\A^\rho\}}: \R^2 \to M_n(\C)$ be defined by 
\[
\Psi_{\{\A^\rho\}}(\rho, \ell) = \mathcal{F}(\A^\rho, \ell),
\]
so that the evaluation map $ev_\Psi:\R^2 \to M_n(\C)$ is 
\[
ev_\Psi(\A, \rho, \ell) = \mathcal{F}(\A^\rho, \ell).
\]
Then, taking $r = 1, N = 2, Q = 2(n-k)^2$, the third condition of Theorem \ref{transversality} is satisfied for any $0 \le k \le n-2$. By Proposition \ref{surjective}, the evaluation map is also transverse to $\mathbf{G}_k$ for any $0 \le k \le n-2$. 

Repeating this for the other two maps, and taking intersections, there then exists a residual subset (hence dense) of $\mathcal{P}^1$ such that for any ${\{\A^\rho\}}$ in the set, all conditions from Proposition \ref{conditions} are satisfied. \end{Proof}

\begin{Proof}[of Lemma \ref{simpledense}]
By Proposition \ref{resid}, and density, we may assume without loss of generality that the family $\{ \A^\rho \}$ in Proposition \ref{simpledense} satisfies the conditions from Proposition \ref{conditions} for each $\rho \in \R$. For each such $\A^\rho$, there is at most one $\ell \in \R$ such that $\nu = i\ell$ is a root of $\det \Delta_{\A^\rho} = 0$, and the root is simple. 

By this assumption, there exist $\ep, L > 0$ such that any root $\lambda(\rho, \nu)$ with $|\Re(\lambda)| < \ep$ is simple. Also, by hyperbolicity at infinity, there are no roots with $|\Re(\lambda)| < \ep$ for $\rho \notin [-L,L]$, choosing $L$ sufficiently large, possibly taking $\ep$ less than the $\ep$ in the statment of Lemma \ref{simpledense}. Then any root with $|\Re(\lambda)| < \ep$ can be parameterized as a $C^1$ function of $\rho$, on a maximal open interval $I \subset \R$ such that $|\lambda| < \ep$. Label the set of such parameterizations $\{\lambda_i(\rho)\}$. Note that there can be no more than countably many such parameterizations. 

Then by Sard's theorem, almost every $\gamma \in (-\ep, \ep)$ is a regular value for every $\Re(\lambda_i(\rho))$. Fix one such $\gamma_0 \in (0, \ep)$. 

Define for $t \in \R$ the operator $S_t: (W^{1,1}_\eta(\R, M_n(\C)) \times M_n(\C) )\to (W^{1,1}_\eta(\R, M_n(\C)) \times M_n(\C) )$ by 
\[
S_t(\A^0) = S_t((A^0, K^0(\cdot))) = (A^0, K^0(\cdot)e^{it(\cdot)}).
\]  
One can check that $\Delta_{S_t(\A^0)}(\nu) = \Delta_{\A^0}(\nu - t), \nu \in \C$, so that $S_t$ shifts all roots of the characteristic equation to the right by an amount $t$. 
Now, let the smooth nonnegative function $\gamma: \R \to \R$ equal 0 outside $[-L+1,L+1]$, equal $\gamma_0$ on $[-L,L]$, and never exceed $\gamma_0$ in between. 
Then the family $\{S_{-\gamma(\rho)}(\A^\rho)\}$ can be seen to satisfy conditions (i) and (ii) of Lemma \ref{simpledense}, and additionally, the roots $\lambda_i(\rho) - \gamma(\rho)$ cross the imaginary axis tranversely. This proves  Lemma \ref{simpledense}. \end{Proof}

\begin{Lemma}\label{forsimple} Let $\{ \A^\rho\} \in \mathcal{P}^1$ be such that $NH(\{\A^\rho\})$ is a finite set, and such that it has only simple crossings. Then 
\[
\iota(\A^+, \A^-) = \iota_c(\A^+, \A^-) = - \textrm{cross}(\{ \A^\rho\}).  
\] 
\end{Lemma}

\begin{Proof} The proof of this lemma follows identically as in \cite{faye2014fredholm}, from analogues of Propositions 4.7 and 4.8, there, and the cocycle property, replacing $L^2$ with $L^\infty$, $C^0$ where necessary; we refer to  \cite{faye2014fredholm} and  \cite{mallet1999fredholm} for proofs. We note that the proof eventually reduces to establishing that the operator $(\frac{d}{d\xi} - i\ell)(\frac{d}{d\xi} + \omega)^{-1}$ is Fredholm index -1 on exponentially decaying spaces $C^0_\gamma(\R, \C^n), L^\infty_\gamma(\R,\C^n)$, for $\gamma < \eta$ and $\omega, \ell > 0$, which can be explicitly verified. \end{Proof}

\begin{Proof}[of Theorem \ref{arho}] From Lemma \ref{simpledense}, we have that, generically, paths cross the axis with only finitely many crossings, all of which are simple. Lemma \ref{forsimple} then gives us that for such a path of operators, the Fredholm index is given by the crossing number.  Putting this all together, we see that Theorem \ref{arho} is proved. \end{Proof} 

\section{Nonlocal Center Manifolds in $C^0$-based Spaces}
\label{s:4}

We  consider the following nonlinear, nonlocal equation: \begin{equation}\label{e:main} -u + K*u + \mathcal{F}(u) = 0 ,\end{equation}
where $K$ is a matrix convolution kernel and $\mathcal{F}(u)(x) = f(u(x)), f \in C^{k}(\mathcal{U}, \R^n)$ a pointwise nonlinearity, $k \ge 1$, for $\mathcal{U}$ a neighborhood of $0 \in \R^n$. We denote $\T u = -u + K*u$. Assuming that $f(0)=0,\,f'(0)=0$, we are interested in small solutions $u(x)$, $\| u \|_{L^\infty}<\delta\ll 1$.  To leading order, one expects that the linearization predicts behavior of small solutions. This fact is commonly captured in center manifold theorems or Lyapunov-Schmidt reduction techniques. For the nonlocal equation \eqref{e:main}, such a reduction was found in  \cite{fs2018center}, parameterizing the set of (possibly weakly) bounded solutions to this equation over the kernel of the linearization. Different from \cite{fs2018center}, we wish to pursue that same goal but relying on $C^0$-based instead of $L^2$-based spaces. We refer to this construction, that we also describe in more detail below, as a center manifold for nonlocal equations. 

As is standard in center manifold constructions, we first use a cut-off function to construct a modified nonlinearity, so that we can use a fixed-point argument in spaces allowing for mild exponential growth. 
%
%
We then show that the set of small bounded solutions to \eqref{e:main} can be described by solutions to a reduced differential equation. This equation is posed on the abstract finite-dimensional vector space given by the kernel of the linearization, allowing for explicit computations of Taylor jets in a straight-forward fashion, using only moments of $K$ and the Taylor series of $f$; see \cite[\S2.6]{fs2018center}. The key to constructing the reduced vector field is this: the analogy of a flow in phase space is the shift operator $u(\cdot) \mapsto u(\cdot + x)$ in function space. This linear shift operator, acting on the nonlinear set of bounded solutions, induces a nonlinear flow when projected onto the kernel. This flow can then be differentiated to obtain a reduced vector field. 

To obtain optimal regularity, we perform a center manifold reduction for the equation with a slightly different nonlinearity,   
\begin{equation}\label{e:v}
  -v + K*v + K*\mathcal{G}(v) = 0,
\end{equation}
with the same assumptions on $\mathcal{G}$ as $\mathcal{F}$. Assuming that $K$ has a derivative, as assumed throughout in Sections \ref{s:2}--\ref{s:3}, we find that small bounded solutions $v\in C^0$ to \eqref{e:v} are automatically small and bounded in $C^1$.  Equations \eqref{e:main} and \eqref{e:v} are equivalent through the change of variables $v = u - f(u)$; starting with \eqref{e:main}, we obtain \eqref{e:v} with  $g(v) = (\mathrm{Id} - f)^{-1}(v) - v$. By the inverse function theorem, $g$ is as smooth as $f$, $g \in C^k(\mathcal{U},\R^n)$. We note, however, that, assuming $f\in C^k$, this $C^k$-change of variables would a priori only yield a $C^{k-1}$ vector field, so that from the perspective of regularity theory, the two formulations may not be equivalent. The formulation \eqref{e:v} yields optimal regularity of the center manifold, while recovering the regularity in \cite{fs2018center} for the $u$ formulation \eqref{e:main}.



\subsection{Hypotheses for Center Manifold Existence}

We require localization of the kernel and smallness of the nonlinearity near the trivial solution: 

\begin{Hypothesis}[Exponentially localized convolution]\label{h:k}
We assume that the matrix convolution operator is exponentially localized and differentiable,  $K \in W^{1,1}_{\eta_0}(\R, M_n(\R))$ for some $\eta_0 > 0$. 
\end{Hypothesis}

\begin{Hypothesis}[Small nonlinearity]\label{h:g}
We asssume that the nonlinearity is small near the origin in the sense that $g\in C^k(\mathcal{U},\R^n)$ for some neighborhood $\mathcal{U}$ of $0\in\R^n$, $1\leq k<\infty$, $g(0)=0$, and $g'(0)=0$.
\end{Hypothesis}

In order to state our main result, we define  the Banach space $C^0_{\sigma}(\R,\R^n)$, for any $\sigma \in \R,$ to be the space $\{ v \in C^0(\R,\R^n) \ | \ \| v(\cdot)e^{\sigma| \cdot |} \|_{C^0} < \infty \}$, and let $C^1_{\sigma}(\R,\R^n)$ be defined analogously. We will often refer to these spaces with $\sigma = -\eta$, simply by $C^0_{-\eta}$ or $C^1_{-\eta}$ for brevity.  

By Hypothesis \ref{h:k}, $\T$ is a bounded operator on $C^0_{-\eta}(\R,\R^n), 0 < \eta < \eta_0$, suppressing notationally the dependence on $\eta$. Moreover, as we will see below, the kernel $\mathcal{E}_0$ of $\T$ is finite-dimensional and independent of $\eta$ for $\eta_0$ sufficiently small, in the sense that the bounded inclusions $\iota_{\eta, \eta'}: C^0_{-\eta} \to C^0_{-\eta'}, \eta < \eta'$, provide kernel isomorphisms.

One can readily see, for example from \cite[\S2.5]{fs2018center}, that there exists a projection operator $\mathcal{Q}: C^0_{-\eta} \to C^0_{-\eta}$
onto the kernel $\mathcal{E}_0$ of $\T$ satisfying $\mathcal{Q}\iota_{\eta, \eta'} = \iota_{\eta,\eta'}\mathcal{Q}$. This projection will play an essential role in our construction of the reduced flow below. 

Next, define the translation operator $\tau_\xi,$ for $\xi \in \R$, by 
\[
(\tau_\xi\cdot v)(x) := v(x-\xi). 
\]
Again, slightly abusing notation, we use the same symbol for the shift on different function spaces. Clearly $
\tau_\xi$ is a bounded operator on $C^0$ and $ C^0_{-\eta}$ for fixed $\xi$. 

We will also use a modified nonlinearity, cutting off $g$ outside a small neighborhood of the origin. Therefore, define $g^\ep:\R^n\to\R^n$ through
\[
        g^\ep(v)=g\left(\chi(\|v\|/\ep)\cdot v\right), 
\]
where $\chi\in C^\infty(\R_{\ge 0},\R)$ is a smoothed version of the indicator function of $[0,1]$,
\[
\chi(t) = \begin{cases} & 1 \textrm{ for }0\leq  t  \le 1 \\ & 0 \textrm{ for }  t  \ge 2
\end{cases}, \qquad\chi(t) \in [0,1].
\]
Denote by $\mathcal{G}$ and $\mathcal{G}^\ep$ the superposition operators associated with $g$ and $g^\ep$, respectively. One readily verifies that the Lipschitz constant of $g^\ep$ is small for $\ep$ small.
Other modifications such as cut-off outside of the nonlinearity or cut-off operators are also allowed as long as the modified nonlinearity possesses a globally small Lipshitz constant.

\subsection{Existence of a Center Manifold}

We are thus ready to state the main center manifold reduction result. In doing so we study solutions to both the unmodified and modified nonlocal equations,
\begin{equation}\label{e:v2}
\T v + K*\mathcal{G}(v) = 0,
\end{equation}
\begin{equation}\label{e:mainmodified}
    \T v + K*\mathcal{G}^\ep(v) = 0.
\end{equation}
\begin{Theorem}\label{CMexistence}
Assume Hypotheses \ref{h:k} and \ref{h:g} on the kernel $K$ and nonlinearity $g$. Recall the definition of the kernel $\mathcal{E}_0$ of $\T$, the projection $Q$ on the kernel, the shift $\tau_\xi$, and the modified nonlinearity $\mathcal{G}^\ep$. 
Consider equations \eqref{e:v2} and \eqref{e:mainmodified}.

Then for all $\eta>0$ sufficiently small, there exist $\ep, \delta >0$, and a map 
\[ 
\Psi: \mathcal{E}_0 \subset C^0_{-\eta}(\R,\R^n) \to \ker\mathcal{Q} \subset C^0_{-\eta}(\R,\R^n),\] with graph \[\mathcal{M} := \{ v_0 + \Psi(v_0) \ | \ v_0 \in \ker\T\} \subset C^0_{-\eta}(\R,\R^n)
,\] 
such that the following hold: 
\begin{enumerate}[(i)]
\item (smoothness and tangency) $\Psi \in C^k,$ with $k$ as in Hypothesis \ref{h:g}, $\Psi(0) = 0$, $D\Psi(0) = 0$;
\item (global center manifold reduction) $\mathcal{M}$ consists precisely of the solutions in $C^0_{-\delta}(\R,\R^n)$ of the modified equation \eqref{e:mainmodified};
\item (local center manifold reduction) any solution $v \in C^0_{-\delta}(\R,\R^n)$ of the unmodified equation \eqref{e:v2} with $\sup_{x \in \R}|v(x)| \le \ep$ is contained in $\mathcal{M}$;
\item (translation invariance) the shift $\tau_\xi$, $\xi \in \R$, acts on $\mathcal{M}$ and induces the reduced flow $\Phi_\xi: \mathcal{E}_0 \to \mathcal{E}_0$ through $\Phi_\xi = \mathcal{Q} \circ \tau_\xi \circ \Psi;$
\item (reduced vector field) the reduced flow $\Phi_\xi(v_0)$ is of class $C^k$ in $v_0, \xi$ and generated by a reduced vector field $h$ of class $C^{k}$ on the finite-dimensional vector space $\mathcal{E}_0.$
\end{enumerate}
In particular, small solutions to $u' = h(u)$ on $\mathcal{E}_0$ are in one-to-one correspondence with small bounded solutions of \eqref{e:v2}.
\end{Theorem}

We refer to the discussion in \cite{fs2018center} and \cite{bakker2018hamiltonian} for further properties of flows on the center manifold, such as dependence on parameters, the computation of Taylor expansions, symmetries and reversibility, Hamiltonian and gradient-like structure, or normal forms. 

We reiterate here that the use of $C^0$-based spaces allows us to obtain optimal regularity of the center manifold and the reduced vector field when compared to the results in \cite{fs2018center}. We also note that the cut-off procedure outlined here is significantly easier than the construction in \cite{faye2020corrigendum} and may well prove more versatile in applications to more complicated, nonlocal nonlinearities. 

\subsection{Proof of Theorem \ref{CMexistence}}
The proof generally follows the strategy in \cite{fs2018center}. We collect properties of the nonlinearity, first, and then study Fredholm properties of the linearization. We then prove existence and regularity of the parameterization of the center manifold using contraction principles on scales of Banach spaces. Lastly, we establish existence and smoothness of the reduced vector field by showing additional smoothness of solutions using bootstraps and then investigating the flow induced by translations of bounded solutions. 

We start by collecting some properties of the superposition operator induced by $g^\ep$: 
\begin{itemize}
         \item $\mathcal{G}^\ep$  is continuous from $C^0_{-\zeta}$ to $C^0_{-\eta}$ for $\zeta, \eta > 0$; moreover, since $g \in C^1(\mathcal{U})$, $\mathcal{G}^\ep$ is Lipschitz in $u$ if $\eta \ge \zeta$, with $\lip_{C^0_{-\zeta} \to C^0_{-\eta}}(\mathcal{G}) \le \lv g^\ep \rv_{C^1} = o_\ep(1)$;
        \item $\mathcal{G}^\ep$  is $k$ times Frechet differentiable from $C^0_{-\zeta}$ to $C^0_{-\eta}$ for $0 < k\zeta < \eta$; 
        \item $\mathcal{G}^\ep(0) = 0$ and, when defined, $D_v\mathcal{G}^\ep(0) = 0$;
        \item $\mathcal{G}^\ep$ is translation-invariant; that is, $\tau_\xi \circ \mathcal{G}^\ep = \mathcal{G}^\ep \circ \tau_\xi$;
\end{itemize}
see for instance \cite{VdBVH1987}. 
We next collect information on the linearization in exponentially weighted spaces. Consider the linear operator \begin{equation*} \mathcal{T}: C^0_{-\eta} \to C^0_{-\eta}, \ \ \ \mathcal{T}(v) = -v + K*v ,\end{equation*}
and its associated characteristic function 
\[
d(\nu) = \det(\mathbb{I}_n + \wh{K(\nu})), \nu \in \C.
\]
The following result determines Fredholm properties in terms of roots of $d$ on the imaginary axis. In fact, the sum of  multiplicities of roots of $d(i\ell)$ on $\ell\in\R$ is  finite. To see this, one first exploits that $\wh{K}$ and thereby $d$ are analytic so that roots have locally finite multiplicities, by exponential localization of $K$. One then notes that  $\wh{K}(i\ell)$ decays as $|\ell|\to\infty$ by regularity of $K$, so that $d$ does not vanish for large $\ell$.  

\begin{Proposition}
Assuming Hypothesis \ref{h:k}, the operator $\mathcal{T}$ is Fredholm with index $M<\infty$, where $M$ is the sum of the multiplicities of roots of $d(i\ell) = \det(\mathbb{I}_n + \wh{K(i\ell})), \ell \in \R$.
\end{Proposition}

\begin{Proof} 
We start by first conjugating $\T$
with the multiplication operator $v(x) \mapsto \cosh(\eta x)\cdot v(x)$ to obtain an operator on $C^0$ of the form considered in Theorem \ref{indexflow}. We note that this theorem refers to operators on complex function spaces, but the corresponding statement for real operators is obtained immediately by restricting to real subspaces. By Theorem \ref{indexflow}, the conjugated operator is Fredholm, with index equal to the number of roots of its characteristic equation that cross the imaginary axis, counted with multiplicity. This quantity is exactly equal to the number of roots $M$ of $\det(I + \wh{K(i\ell}))$, counted with multiplicity, for $\ell \in \R$, so the proposition follows. \end{Proof}

We now define the bordered operator 
\begin{equation*} 
  \wt{\T}: C^0_{-\eta} \to C^0_{-\eta} \times \mathcal{E}_0,\qquad \wt{\T}(v) = (\T(v),\mathcal{Q}(v))
,\end{equation*}
which, when solving $ \wt{\T}v=(f,v_0)$, forces  $\mathcal{Q}v = v_0$ for a given $v_0 \in \mathcal{E}_0$, the kernel of $\T$.

Fredholm bordering theory guarantees that $\wt{\T}$ is Fredholm, since a finite number of dimensions are being added onto the range, and has index 0, since $M = \dim(\ker\T)$. Furthermore, it is now one-to-one, since $\T v = 0$ and $\mathcal{Q}v = 0$ imply $v = 0$. Therefore the bordered operator is in fact invertible with bounded inverse, such that 
\begin{equation*}\lv\wt{T}^{-1}\rv_{\mathcal{L}(C^0_{-\eta},C^0_{-\eta})} \le C(\eta),\end{equation*}
for a constant $C(\eta)$, with $C(\eta)$ continuous in $\eta$ for $0 < \eta < \eta_0$.

We are now able to set up a fixed point equation using the bordered equation
\begin{equation}\label{e:bordered}
    \wt{\T}(v) + \wt{\mathcal{G}}^\ep (v,v_0) = 0,
\end{equation}
 where $\wt{\mathcal{G}}^\ep(v,v_0) = (\mathcal{G}^\ep(v),-v_0)$; note that this is equivalent to the original equation. Rewriting \eqref{e:bordered}, we find for a given $v_0 \in \mathcal{E}_0$,
\begin{equation}
    v = -\wt{\T}^{-1}(\wt{\mathcal{G}}(v,v_0)).
\end{equation}
We view this equation as a fixed-point equation on $C^0_{-\eta}$ with parameter $v_0$. We claim that the map $-\wt{\T}^{-1}(\wt{\mathcal{G}}(\cdot,v_0))$ is a contraction mapping. To see this, we use that $g(0)=0$ and $g'(0)=0$ to find that 
\begin{align*}
    \delta_0(\ep) &:= \sup_{v \in C^0_{-\eta}} \lv \mathcal{G}^\ep(v) \rv_{C^0_{-\eta}} = o(\ep)\\
    \delta_1(\ep) &:= \textrm{Lip}_{C^0_{-\eta}\to C^0_{-\eta}}(\mathcal{G}^\ep ) = o_{\ep}(1),
\end{align*}


which in turn implies that 
\begin{align*}
   \lv \wt{\T}^{-1}(\wt{\mathcal{G}}(v,v_0)) \rv_{C^0_{-\eta}}&\le C(\eta)\left(\delta_0(\ep) + \lv v_0 \rv_{C^0_{-\eta}}\right) \\
\lv \wt{\T}^{-1}(\wt{\mathcal{G}}(v_1,v_0)) - \wt{\T}^{-1}(\wt{\mathcal{G}}(v_2,v_0)) \rv_{C^0_{-\eta}} &\le C(\eta)\delta_1(\ep)\lv v_1 - v_2 \rv_{C^0_{-\eta}}  
\end{align*}
for all $v, v_1, v_2 \in C^0_{-\eta}$ and $v_0 \in \mathcal{E}_0$. 

Then, letting $\overline{\eta} \in (0,\eta_0)$ and $\wt{\eta} \in (0, \frac{\overline{\eta}}{k})$, for $\ep$ sufficiently small, we have $C(\eta)\delta_1(\ep)<1$ for $\eta \in [\wt{\eta},\overline{\eta}]$, so that $\wt{\T}^{-1}(\wt{\mathcal{G}}(\cdot,v_0))$ defines a contraction mapping on $C^0_{-\eta}$, and has a unique fixed point $v = \Phi(v_0)$. Since the fixed point iteration is Lipschitz in $v_0$, the map $\Phi$ is also Lipschitz, with $\Phi(0) = 0$ because the fixed point is unique. For each $\eta,$ this then defines a Lipshitz map $\Psi: \mathcal{E}_0 \to \ker \mathcal{Q}$ such that 
\begin{equation*}
    \Phi(v_0) = v_0 + \Psi(v_0).
\end{equation*}
Note that $\Phi$ commutes with translations $\tau_\xi$ by uniqueness of the fixed point. 

We next turn to smoothness of $\Phi$, following ideas in  \cite{VdBVH1987}.
\begin{Proposition}\label{Ck_of_Phi}
Under the same assumptions as Theorem \ref{CMexistence}, for each $1 \le p \le k$ and for each $\eta \in (p\wt{\eta},\overline{\eta})$, the map $\Psi$ is $C^p$ from $\mathcal{E}_0$ to $C^0_{-\eta}$.
\end{Proposition}

In order to prove this, we recall the following result from \cite{VdBVH1987} on contractions on scales of embedded Banach spaces.

Let $\mathcal{X}, \mathcal{Y}, \mathcal{Z}$ and $\Lambda$ be Banach spaces with norms denoted by $\lv \cdot \rv_{\mathcal{X}}, \lv \cdot \rv_{\mathcal{Y}}, \lv \cdot \rv_{\mathcal{Z}}$ and $\lv \cdot \rv_{\Lambda}$, with continuous embeddings 
\[
\mathcal{X} \overset{\mathcal{J}}{\hookrightarrow} \mathcal{Y} \overset{\mathcal{I}}{\hookrightarrow} \mathcal{Z}.
\]
Consider the fixed point equation \begin{equation}\label{e:fixedpt} y = \textbf{f}(y, \lambda), \end{equation}
where $\textbf{f}: \mathcal{Y} \times \Lambda \to \mathcal{Y}$ satisfies the following conditions: 
\begin{enumerate}[(i)]
    \item $\mathcal{I}\textbf{f}: \mathcal{Y} \times \Lambda \to \mathcal{Z}$ has continuous partial derivative $D_y(\mathcal{I}\textbf{f}): \mathcal{Y} \times \Lambda \to \mathcal{L}(\mathcal{Y}, \mathcal{Z})$ with \[D_y(\mathcal{I}\textbf{f})(y,\lambda) = \mathcal{I}\textbf{f}^{(1)}(y,\lambda) = \textbf{f}_1^{(1)}(y,\lambda)\mathcal{I}, \ \ \ \textrm{ for all } (y,\lambda) \in \mathcal{Y} \times \Lambda, \] for some $\textbf{f}^{(1)}: \mathcal{Y} \times \Lambda \to \mathcal{L}(\mathcal{Y})$ and $\textbf{f}_1^{(1)}: \mathcal{Y} \times \Lambda \to \mathcal{L}(\mathcal{Z})$. 
    \item $\textbf{f}_0: \mathcal{X} \times \Lambda \to \mathcal{Y}, (y_0, \lambda) \mapsto \textbf{f}_0(y_0, \lambda) = \textbf{f}(\mathcal{J}y_0, \lambda)$ has continuous partial derivative $D_\lambda \textbf{f}_0: \mathcal{X} \times \Lambda \to \mathcal{L}(\Lambda, \mathcal{Y}).$
    \item There exists $\kappa \in [0,1)$ such that \[ \lv \textbf{f}(y, \lambda) - \textbf{f}(\wt{y},\lambda)\rv_\mathcal{Y} \le \kappa \lv y - \wt{y}\rv_{\mathcal{Y}}, \ \ \ \textrm{ for all } y, \wt{y} \in \mathcal{Y}, \ \ \textrm{ for all } \lambda \in \Lambda,\]
    and \[ \lv \textbf{f}^{(1)} (y, \lambda) \rv_{\mathcal{L}(\mathcal{Y})} \le \kappa, \ \ \lv \textbf{f}^{(1)}_1(y, \lambda) \rv_{\mathcal{L}(\mathcal{Z})} \le \kappa, \ \ \ \textrm{ for all } (y, \lambda) \in \mathcal{Y} \times \Lambda.\]
    \item Let $y = \wt{y}(\lambda) \in \mathcal{Y}$ be the unique solution of \eqref{e:fixedpt} for $\lambda \in \Lambda$. Suppose that $\wt{y}(\lambda) = \mathcal{J}\wt{y}_0(\lambda)$ for some continuous $\wt{y_0}: \Lambda \to \mathcal{X}$.
\end{enumerate}

These conditions allow consideration of the following equation in $\mathcal{L}(\Lambda, \mathcal{Y}):$
\begin{equation}\label{e:Theta}
    \Theta = \textbf{f}^{(1)}(\wt{y}(\lambda), \lambda)\Theta + D_\lambda \textbf{f}_0(\wt{y}_0(\lambda), \lambda),
\end{equation}
which has a unique solution $\wt{\Theta}(\lambda) \in \mathcal{L}(\Lambda, \mathcal{Y})$ for any $\lambda \in \Lambda$ from condition (iii). The following result is proved in \cite{VdBVH1987}:

\begin{Theorem}\label{scales} Assume conditions (i)-(iv). Then the solution map $\wt{y}: \Lambda \to \mathcal{Y}$ of \eqref{e:Theta} is Lipschitz continuous, and $\mathcal{I}\wt{y}: \Lambda \to \mathcal{Z}$ is of class $C^1$, with \begin{equation} D_\lambda \mathcal{I}\wt{y}(\lambda) = \mathcal{I}\wt{\Theta}(\lambda), \ \ \ \textrm{ for all } \lambda \in \Lambda. \end{equation}
\end{Theorem}

We turn now to the proof of Proposition \ref{Ck_of_Phi}. 
\begin{Proof}[(of Proposition \ref{Ck_of_Phi})]
This argument is a straightforward analogue of the proof of Lemma 6 from \cite{VdBVH1987}, as well as appendix A from \cite{fs2018center}. We begin by letting $p =1$ and fixing $\eta \in (\wt{\eta}, \overline{\eta}].$ Then apply Theorem \ref{scales} with $\mathcal{X} = \mathcal{Y} = C^0_{-\wt{\eta}}, \mathcal{Z} = C^0_{-\eta}, \Lambda = \mathcal{E}_0$ and $\textbf{f}(y, \lambda) = -\wt{\T}^{-1}(\mathcal{G}^\ep(y;\lambda))$. One can check that assumptions (i)-(iv) are verified, so that $\Phi: \mathcal{E}_0 \to C^0_{-\eta}$ is of class $C^1$ with derivative $\Phi^{(1)}(v_0) := D\Phi(v_0) \in \mathcal{L}(\mathcal{E}_0, C^0_{-\eta})$ the unique solution of \begin{equation}\label{e:firstderiv} \Theta = D_y\textbf{f}(\Phi(v_0), v_0)\Theta + D_\lambda\textbf{f}(\Phi(v_0), v_0) := F_1(\Theta, v_0).\end{equation}
Now, the mapping $F_1: \mathcal{L}(\mathcal{E}_0,C^0_{-\eta}) \times \mathcal{E}_0 \to \mathcal{L}(\mathcal{E}_0, C^0_{-\eta})$ is a uniform contraction for each $\eta \in [\wt{\eta},\overline{\eta}]$, 
so the fixed point of \eqref{e:firstderiv} belongs in fact to $\mathcal{L}(\mathcal{E}_0,C^0_{-\wt{\eta}})$.  The mapping $\Phi^{(1)}: \mathcal{E}_0 \to \mathcal{L}(\mathcal{E}_0, C^0_{-\eta})$ is continuous if $\eta \in (\wt{\eta},\overline{\eta}].$

If $k \ge 2$, we now continue by induction. Let $1 \le p < k$, and suppose that for all $q$ with $1 \le q \le p$ and for all $\eta \in (q\wt{\eta},\overline{\eta}]$ the mapping $\Phi: \mathcal{E}_0 \to C^0_{-\eta}$ is of class $C^p$, with $\Phi^{(q)}(u_0) := D^q\Phi(v_0) \in \mathcal{L}^{(q)}(\mathcal{E}_0, C^0_{-q\wt{\eta}})$ for each $v_0 \in \mathcal{E}_0$ and $\Phi^{(q)}: \mathcal{E}_0 \to \mathcal{L}^{(q)}(\mathcal{E}_0, C^0_{-\eta})$ continuous if $\eta \in (q\wt{\eta}, \overline{\eta}]$. Suppose in addition that $\Phi^{(p)}(v_0)$ is the unique solution of an equation that is of the form \begin{equation}
    \Theta^{(p)} = D_y\textbf{f}(\Phi(v_0),v_0)\Theta^{(p)} + H_p(v_0) := F_p(\Theta^{(p)}, v_0),
\end{equation}
with $H_1(u_0) = D_\lambda(\textbf{f}(\Phi(u_0),u_0)$ and, for $p \ge 2,H_p(u_0)$ is given as a finite sum of terms of the form \[D_y^{(q)}\textbf{f}(\Phi(v_0),v_0)(D^{r_1}\Phi(v_0), ..., D^{r_q}\Phi(v_0)), \] with $2 \le q \le p, 1 \le r_i < p$ for all $i = 1, ... , q,$ and $r_1 + ... + r_q = p$. By similar reasoning as before, we note that $H_p(u_0) \in \mathcal{L}^{(p)}(\mathcal{E}_0, C^0_{-p\wt{\eta}}).$ Therefore, $F_p:\mathcal{L}^{(p)}(\mathcal{E}_0, C^0_{-p\wt{\eta}}) \times \mathcal{E}_0 \to \mathcal{L}^{(p)}(\mathcal{E}_0, C^0_{-p\wt{\eta}})$ is well defined and a uniform contraction for $\eta \in [p\wt{\eta},\overline{\eta}].$ However, the term $D_y\textbf{f}(\Phi(v_0))\Phi^{(p)}$ is not continuously differentiable, either with respect to $\Phi^{(p)}$ or the parameter $u_0$, so we apply Theorem \ref{scales} with three different Banach spaces. Let $\eta \in ((p+1)\wt{\eta},\overline{\eta}]$,$\sigma \in (\wt{\eta},\frac{\eta}{(p+1)})$, and $\zeta \in ((p+1)\sigma, \eta).$ We need to show that the hypotheses of Theorem \ref{scales} are satisfied with $\mathcal{X} = \mathcal{L}^{(p)}(\mathcal{E}_0, C^0_{-p\sigma}), \mathcal{Y} = \mathcal{L}^{(p)}(\mathcal{E}_0, C^0_{-{\zeta}})$, and $\mathcal{Z} = \mathcal{L}^{(p)}(\mathcal{E}_0, C^0_{-{\eta}})$, $\Lambda = \mathcal{E}_0$ and $\textbf{f} = F_p$. Condition (iii) holds because $C(\eta)\delta_1(\ep) < 1$ for $\eta \in [\wt{\eta},\overline{\eta}]$. Condition (iv) holds by the induction hypothesis and because $\sigma > \wt{\eta}$. Now, the map $D_y\textbf{f}(\Phi(v_0),v_0)$ is continuous from $\mathcal{E}_0$ into $\mathcal{L}^(C^0_{-\zeta}, C^0_{-{\eta}})$, because $\Phi: \mathcal{E}_0$ is continuous and $\eta > \zeta$ (see \cite{VdBVH1987}, Lemma 4). Further, by the same, $D_y\textbf{f}(\Phi(v_0),v_0)$ is $C^1$ from $\mathcal{E}_0$ into $\mathcal{L}(C^0_{-p\sigma},C^0_{-\zeta})$, because $\zeta > (p+1)\sigma$ and $\Phi \in C^1$. It thus remains to show that $H_p: \mathcal{E}_0 \to C^0_{-\zeta}$ is of class $C^1$. This again follows by the same reasoning as \cite{VdBVH1987}, Lemma 7. Then we can use Theorem \ref{scales} and conclude that $\Phi^{(p)}: \mathcal{E}_0 \to \mathcal{L}^{(p)}(\mathcal{E}_0, C^0_{-{\eta}})$ is of class $C^1$ and hence $\Phi: \mathcal{E}_0 \to C^0_{-\eta}$ is of class $C^{p+1}$ if $\eta \in ((p+1)\wt{\eta},\overline{\eta}].$
\end{Proof}

\subsection{Existence of a reduced vector field}

The next step in the proof of Theorem \ref{CMexistence} is the construction of the reduced vector field. As mentioned in the introduction, this is obtained by differentiating the action of the shift operator projected onto the kernel. Therefore, to start with, we would like to show that the solutions in the center manifold in fact belong to $C^1_{-\eta}$, so that the shift map can be differentiated. 


For $v \in C^0_{-\eta}$ a solution of \eqref{e:v2}, we have
\begin{equation}\label{e:bootstrap} 
   v(x) = (K*(\mathrm{Id}+\mathcal{G}^\ep)(v))(x).
\end{equation}
Now, the map $(\mathrm{Id} + \mathcal{G}^\ep)$ is, as a superposition operator, a $C^k$ map from $C^0_{-\zeta}$ to $C^0_{-\eta}$, as proved in \cite{VdBVH1987}, for $0 < k\zeta < \eta$. The map $u \mapsto K*u$ is a bounded linear map from $C^0_{-\eta}$ to $C^1_{-\eta}$, due to the fact that $K \in W^{1,1}_{\eta}$. Then we have that $v \in C^1_{-\eta}$, with the composition $K*((\mathrm{Id} + \mathcal{G}^\ep)\circ(\mathrm{Id}+ \Psi))$ a $C^k$ map from $\mathcal{E}_0$ to $C^1_{-\eta}$. 
%
%
%
%
%
%
%
%
%
%

Now, consider the action of the shift operator 

\begin{align*}
\R \times C^1_{-\eta} &\to C^0_{-\eta} \\
(x,u) &\mapsto \tau_x v = v(\cdot + x). \end{align*}

 We have that for a given $x$, $\tau_x$ is a bounded linear operator which maps bounded solutions of \eqref{e:main} to bounded solutions. The following commutative diagram shows how $\tau_x$ induces a flow on the kernel: 


 \begin{center}
$\xymatrix @!=2.5cm{
\mathcal{E}_0 \ar[r]^{\mathrm{Id}+\Psi} \ar@{}[drrr]|*{\fontsize{15cm}{15.5cm}\selectfont \circlearrowright} \ar[d]_*{\varphi_x} & C^0_{-\zeta}\ar[r]^{\mathrm{Id} + \mathcal{G}^\ep} & C^0_{-\eta} \ar[r]^{K*}& C^1_{-\eta}\ar[d]^*{\tau_x} \\
\mathcal{E}_0 \ar@<3pt>[rrr]^{\iota_{\zeta,\eta} \circ (\mathrm{Id}+\Psi)} & & &  C^0_{-\eta} \ar@<3pt>[lll]^{{\mathcal{Q}}}}
 $  
 \end{center}

The diagram commutes because the composition $K*(\mathrm{Id} + \mathcal{G}^\ep)$ is the identity on the image of $\mathrm{Id} + \Psi$. Now, $\tau_x$ is bounded linear, as well as continuously differentiable in $x$, with derivative equal to the bounded linear map $v(\cdot + x) \mapsto v'(\cdot + x)$. Then the composition $\mathcal{Q}\circ \tau_x \circ K*((\mathrm{Id} + \mathcal{G}^\ep)\circ(\mathrm{Id} + \Psi))(\cdot)$ is also continuously differentiable in $x$, since $\mathcal{Q}$ is a bounded linear projection. The maps $\Phi$ and $(\mathrm{Id} + \mathcal{G}^\ep)$ are each $C^k$ on their respective function spaces, so that $\varphi_x$ inherits the regularity of the composition, and $\frac{d\varphi_x}{dx}|_{x=0}$ is thus a $C^k$ vector field on $\mathcal{E}_0,$
\begin{equation}\label{e:CMreduced}
    \frac{d\varphi_x}{dx}|_{x=0}:= h(x).
\end{equation}
Likewise, solutions to $\frac{dv}{dx} = h(x), v(0) = v_0$ yield trajectories $\varphi_x(v_0)$ and solutions  $(\mathrm{Id} + \Phi)(\varphi_x(v_0))$ to the nonlocal equation. 

Thus small bounded solutions to \eqref{e:main} can be obtained through solutions to a reduced differential equation on the finite-dimensional kernel, which is in turn obtained by differentiating the reduced flow at $x = 0$.

\subsection{Reduced Vector Field in Original Coordinates} 

The reduced vector field corresponding to the original coordinates can be found by repeating the above procedure with the map $(\mathrm{Id} + g) \circ K*$ instead of just $K*$, and $\iota_{\sigma,\eta} \circ (\mathrm{Id} + g) \circ(\mathrm{Id} + \Psi)$ in the place of $\iota_{\zeta,\eta} \circ (\mathrm{Id} + \Psi)$. In other words, the shift action on the $u$- rather than the $v$-coordinates is differentiated. This will yield a $C^{k-1}$ vector field, since the change-of-coordinate map $(\mathrm{Id} + g)$ is only $C^{k-1}$ from $C^1_{-\eta}$ to $C^1_{-\sigma}$, $0 < (k+1)\eta < \sigma$. This nevertheless recovers the smoothness of the reduced vector field in \cite{fs2018center}, since their $C^k$ reduced vector field corresponded to a $C^{k+1}$ pointwise nonlinearity. 








\section{Application of Center Manifolds: a $C^1$ Lyapunov-Center Theorem}

As an application to Theorem \ref{CMexistence}, we consider the following equation: 
\begin{equation}\label{e:nonlin} 0 = -u + k*(Au + N(u)), \end{equation}
where $A \in GL_n(\R), u \in C^0(\R, \R^n), k \in W^{1,1}(\R, M_n(\R))$, with $k(-x) = k(x)$, and $N(u)(x) = f(u(x))$ a pointwise nonlinearity given by $f \in C^1(\mathcal{V}, \R^n)$, $\mathcal{V}$ a neighborhood of $0 \in \R^n$, with $f(0) = f'(0) = 0$.

The assumption that $k$ be even is intended to be reminiscent of a reversibility condition for nonlinear ODEs. In this context, Lyapunov-Center theorems are a well-known set of results for reversible systems. In essence, they say that if the linearized problem at a given equilibrium has purely imaginary eigenvalues $i\omega_*$ that are non-resonant in a certain sense, then there exists a family of periodic solutions nearby for the full nonlinear problem. Moreover, this family is parameterized roughly by the positive real amplitude and shift parameters. We seek here to establish such a theorem in a nonlocal, spatial dynamics setting, where eigenvalues now correspond to roots of $d(\nu) = \det(\Id +\wh{k}(\nu)A)$. Our main emphasis is on proving that the family of periodic solutions comprises \emph{all} small bounded solutions when $\pm i \omega_*$ are the only roots on the imaginary axis and simple, with minimal assumptions on the regularity of the nonlinearity.

\begin{Hypothesis}\label{h:om}
Assume that there exists $\omega_* > 0$ such that $d(i\omega_*) = \det(\Id + \wh{k}(i\omega_*)A) = 0$, and that $d'(i\omega_*) \neq 0$. Additionally assume that $d(i\omega) \neq 0$ for $\omega \notin \omega_*\Z$. 

\end{Hypothesis}

\begin{Theorem}\label{per_existence} Assuming Hypothesis \ref{h:om}, there exists $\delta > 0$, a continuous frequency function $\omega:[0,\delta)\to\R$ with   $\omega(0) = \omega_*$,
and a 2-dimensional family of periodic solutions to \eqref{e:nonlin}, 
\begin{align*}
u_c: [0,\delta)\times [0,2\pi) &\to \  C^0(\R,\R^n) \\ 
(a,\tau) \qquad  &\mapsto \  u_c\left(\omega(a)(\cdot + \tau);a\right),
\end{align*}
with 
\[
a\mapsto u_c(\cdot;a)\in C^0(\R,\R^n) \text{ continuous},\quad   
u_c(y + 2\pi; a) = u_c(y;a) \textrm{, and } u_c(y;0) = 0.
\]
\end{Theorem}

This theorem, combined with Theorem \ref{CMexistence}, will allow us to prove the following: 

\begin{Theorem}[Nonlocal Lyapunov-Center Theorem]\label{LC} Assume the conditions of Hypothesis 5.1. Assume also that $d(i\omega) \neq 0$ for $|\omega| \neq \omega_*$, and that $k \in W^{1,1}_{\eta_0}$ for $\eta_0 > 0$. Then there exists $\ep > 0$ such that all solutions $u$ to \eqref{e:nonlinper} satisfying $\lv u \rv_{C^0} < \ep$ are periodic and given by the family found in Theorem \ref{per_existence}. 
\end{Theorem}

\begin{Remark}[Necessity of linear conditions]
It is well known that resonances can destroy families of periodic orbits with frequencies that possess higher harmonics. On the other hand, the presence of other roots gives non-uniqueness of periodic families already in the linear case. Lastly, the presence of multiple roots, $d'(i \omega_*)=0$ usually leads to existence of invariant tori, heteroclinic, and homoclinic orbits; see for instance \cite{ioossperoueme} on the reversible Hamiltonian Hopf bifurcation. From this perspective, the assumptions of Theorems \ref{per_existence} and \ref{LC} are necsessary, even for ODEs. 
\end{Remark}

\begin{Remark}[Coherent structures and group velocities]\label{r:LC}
In many contexts, the vanishing of $d'(i\omega_*)$ can be associated with a vanishing group velocity. Consider for example the Kawahara equation in a frame with speed $c>0$,  
\[
u_t=(-\alpha u_{xxxx}+u_{xx}+cu - u^2)_x,
\]
with dispersion relation for solutions $u(t,x)=e^{i(kx-\Omega t)}$ of the linearized equation,
\[
\Omega=\alpha k^5+k^3-ck.
\]
Studying periodic wave trains that are stationary in this frame, we look at $-\alpha u_{xxxx}+u_{xx}+cu - u^2=0$, with characteristic equation $d(i\omega)=-\alpha \omega^4-\omega^2+c$. A root $d(i \omega_*)=0$ gives a root of the dispersion relation with $\Omega=0$ and $k=\omega_*$. The group velocity, $d\Omega/dk$ at this root now vanishes precisely when $d'(i\omega_*)=0$. 

From this perspective, our main result establishes existence of small-amplitude traveling waves as predicted by the linearization, and the absence of any other, possibly non-periodic waves, as long as the group velocity does not vanish in the chosen coordinate frame. We show this absence of non-periodic small traveling waves, such as solitary waves, for minimal assumptions on the regularity of the nonlinearity, noting that continuous differentiability is necessary to give sufficient meaning to the linearization at the origin. In the case when group velocities vanish, existence of nonperiodic waves has been established in many situations, including for instance reductions to KdV or NLS type modulation equations. 
\end{Remark}
We start the remainder of this section with the proof of Theorem \ref{per_existence}, which is essentially proved in four steps: 

\textit{Step 1:} Reduce \eqref{e:nonlin} to a 1-dimensional equation using Lyapunov-Schmidt reduction;

\textit{Step 2:} Set up a contraction argument for the reduced equation;

\textit{Step 3:} Prove contraction properties, yielding a 1-parameter family of solutions;

\textit{Step 4:} Extend the resulting 1-parameter family of solutions to a 2-parameter family by adding a shift parameter.

Theorem \ref{LC} will then follow almost immediately using the center manifold theorem. 

The difficulty in steps 2 and 3 lies in the fact that for a $C^1$ nonlinearity, the reduced equation cannot be solved with the Implicit Function Theorem, since the linear terms vanish. Dividing by the parameter to eliminate the trivial solution does produce linear terms, but loses regularity, so that a more hands-on contraction argument rather than an implicit function theorem is needed to establish existence and uniqueness, taking into account different smoothness in variables and parameters. 

\subsection{Lyapunov-Schmidt Reduction and Derivation of the Reduced Equation}

Let $\wt{u}(x) = u(\omega x)$. Then, changing variables, equation \eqref{e:nonlin} phrased in terms of $\wt{u}$ becomes:
\begin{equation}\label{e:nonlinper}
   0 = -\wt{u} + k_\omega * (A\wt{u} + N(\wt{u})) ,
\end{equation}
where $k_\omega(\cdot) = \frac{1}{\omega} k (\frac{1}{\omega} \cdot)$. We let $F(\omega, u) =-u + k_\omega*(Au + N(u))$, and consider $F(\omega, \cdot)$ as an operator on $C^0_{2\pi, even}(\R,\R^n)$, the set of $C^0$ functions that are $2\pi$-periodic and even. Note that $F$ is a well-defined operator from this function space into itself. 

\begin{Proposition}
 The linearization $\mathcal{L}_{\omega_*} := D_uF(\omega_*, 0)$ is Fredholm index 0, with a 1-dimensional kernel. 
\end{Proposition} 

\begin{Proof} We have that $D_uF(\omega_*, 0)v = -v + k_{\omega_*}*Av.$ The operator $u \mapsto k_{\omega_*}*u$ is compact, because it maps into $C^1_{2\pi, even}(\R, \R^n)$, which is compactly embedded in $C^0_{2\pi, even}(\R, \R^n)$. Then because $D_uF(\omega_*,0)$ is the sum of the identity operator and a compact operator, it is Fredholm, with Fredholm index 0.

We turn now to the kernel. By a calculation, in the space of Fourier series, the linearization of $F$ at $(\omega_*,0)$ is given by $\wh{\mathcal{L}}\wh{u}(j) = (-\Id + \wh{k}(i\omega_*j)A)\wh{u(j)}$. Since we have $d(i\omega_*) = 0$, $ d'(i\omega_*) \neq 0$, then for $j = 1$, the matrix $(\Id + \wh{k(i\omega_*)})$ has a 1-dimensional kernel, spanned by a vector $v_*$, $| v_* |= 1$. For $j \neq 1$, the operator is invertible, since $d(i\omega_*j) \neq 0, j \neq 1$, by assumption. Therefore the kernel of $\mathcal{L}_{\omega_*}$ can be parameterized as $\{ a v_* \cos(x) \ | \ a \in \R \},$ and is 1-dimensional. \end{Proof}

Note that since $\mathcal{L}_{\omega_*}$ is Fredholm index 0, and since $\det(-\Id + A^T\wh{k(i\omega_*)})^T = \det (-\Id + \wh{k(i\omega_*)}A)^T = \det(-\Id + \wh{k(i\omega_*)}) = 0$, then there also exists a vector $v_{ad} \in \R^n, | v_{ad} | = 1$, such that  $\ker\mathcal{L}_{\omega_*}^* = \{ c v_{ad}\cos(x) \ | \ c \in \R\}$.

Now, let $u = a v_* \cos(x) + u_1(x)$, where $u_1(\cdot) \in (\ker\mathcal{L}_{\omega_*})^\perp$, and let $\mathcal{P}$ be the $L^2$-orthogonal projection onto the range of $\mathcal{L}_{\omega_*}$ defined by 
\[
\mathcal{P}u = u - \frac{1}{\pi}\left\langle u(x), v_{ad}\cos(x)\right\rangle_{L^2([0,2\pi], \R^n)} v_{ad}\cos(x).
\] 
Then let 
\begin{equation}\label{e:LSsys}
\begin{split}
F_1(\omega, a, u_1) &= \mathcal{P} F\left(\omega, av_*\cos(x) + u_1(x)\right)  \\
F_0(\omega, a, u_1) &= (1-\mathcal{P})F\left(\omega, av_*\cos(x) + u_1(x)\right),
\end{split}
\end{equation}
with $F_1: \R \times \R \times (\ker(\mathcal{L}_{\omega_*}))^\perp \to \textrm{Ran}(\mathcal{L}_{\omega_*}), F_0:  \R \times \R \times (\ker(\mathcal{L}_{\omega_*}))^\perp \to  \coker(\mathcal{L}_{\omega_*}).$ Since the cokernel of $\mathcal{L}_{\omega_*}$ is one-dimensional, we let $\mathcal{P}_s = \frac{1}{\pi}\left\langle \cdot, v_{ad}\cos(x)\right\rangle_{L^2([0,2\pi]}$, so that $\mathcal{P}_s F_0$ is scalar. Note that $\mathcal{P}_s$ is an isomorphism from $\coker(\mathcal{L}_{\omega_*})$ to $\R$, with $\mathcal{P}_s(1 - \mathcal{P}) = \mathcal{P}_s$. Solutions to the system 
\begin{equation}
\begin{split}
 0 &= F_1(\omega,a,u_1) \\ 0 &= \mathcal{P}_s F_0(\omega,a,u_1) 
\end{split}
\end{equation}
are thus equivalent to solutions to $F(\omega,u) = 0$.

We now exploit Fredholm properties of the linearization to solve $F_1$ near the trivial solution:

\begin{Proposition}
There exists a neighborhood $U$ of $(\omega_0, 0)$ and a $C^1$ function $\psi: U \times \R \to \textrm{Ran}(\mathcal{L}_{\omega_*})$, such that $u_1 = \psi(\omega,a)$ is the unique solution to $F_1(\omega, a, \cdot ) = 0$. Moreover, we have  $\psi(\omega,0) = 0$ and a neighborhood of $(\omega_*,0)$ and a constant $C_1$ such that on that neighborhood, $\|\partial_\omega \psi(\omega,a)\|_{C^0} \le C_1 |a|$. 
\end{Proposition}

\begin{Proof}
We will use the Implicit Function Theorem, for which we will need to establish that $F_1$ \ is \ $C^1$ with respect to $\omega, a, u_1$, and $D_{u_1}F_1(\omega_*, 0, 0)$ is bounded invertible.

First, to show that $F_1$ is $C^1$, we know that $F_1$ is $C^1$ in $a$ and $u_1$ since $N$ is $C^1$, and the remaining terms are linear. As to differentiability in $\omega$, first note that since $k \in W^{1,1}(\R, M_n(\R))$, $k$ is absolutely continuous. One can also check that $\int_{\omega_1}^{\omega_2}  \partial_\omega (k_\omega*u) d\omega < \infty$ for any $\omega_1, \omega_2 > 0$, since $\lv \partial_{\omega}k_\omega \rv_{L^1} = \frac{1}{\omega} \lv k - k' \rv_{L^1}$. Then we will have that $F_1$ is differentiable with respect to $\omega$, with $\partial_\omega F_1u = (\partial_\omega k_\omega)*(Au + N(u)).$ To show that $\partial_\omega F_1$ is continuous in $\omega,$ it suffices to show that the function $\partial_\omega k_\omega = \frac{1}{\omega^2}(k'(\frac{1}{\omega} \cdot) - k(\frac{1}{\omega}\cdot))$ is continuous in $L^1$ in $\omega$. This can be done by finding a compact interval outside of which the tails of $k$ and $k'$ are small enough, and then approximating $k$ and $k'$ inside sufficiently well by continuous functions. 
Since all three partial derivatives are continuous, the function $F_1$ is jointly $C^1$ with respect to $\omega, a, u_1$. As to $D_{u_1}F_1(\omega_*,0,0),$ this is the restriction of $D_{u_1}F(\omega_*,0,0)$, which is Fredholm index 0, to the complement of its kernel, projected onto its range. It will thus be both one-to-one and onto, hence bounded invertible.

As a consequence, by the Implicit Function Theorem, there exists a neighborhood $U$ of $(\omega_0, 0)$ and a $C^1$ function $\psi: U \times \R \to \textrm{Ran}(\mathcal{L}_{\omega_*})$ uniquely solving $F_1(\omega, a, \psi(\omega,a)) = 0$.

It remains to establish the properties of $\psi$ stated. 
The first property is true because $u = 0$ solves the original equation, and because $\psi$ is unique. To justify the second property, by differentiating the equation $F_1(\omega,a,u_1) = 0$ with respect to $\omega$ and using the chain rule, we obtain 
 \[
 \partial_\omega \psi(\omega,a) =  \partial_{u_1}F_1(\omega,a,\psi(\omega,a))^{-1}\partial_\omega F_1(\omega,a,\psi(\omega,a)),
 \]
 provided the inverse exists. However, because the set of invertible linear maps is open, the inverse will exist on some neighborhood $U_1$ of $(\omega_*,0);$ moreover, there exists a uniform bound for a closed subset of that neighborhood.  Now, note that the function \[\partial_\omega F_1(\omega,a,u_1) = \mathcal{P}(-a\cos(x) + u_1(x)) + (\partial_\omega k_\omega)*(A(-a\cos(x) + u_1(x)) + N(-a\cos(x) + u_1(x)),\] while no longer $C^1$ in $\omega$, is still $C^1$ in $a$ and $u_1$, with $\partial_\omega F_1(\omega,0,\psi(\omega,0)) = 0$. Hence we can write $|\partial_\omega F_1(\omega,a,u_1)| \le a \sup_{a,\omega}|\partial_a\partial_\omega F_1(\omega,a,u_1)|$, and since $\partial_a\partial_\omega F_1(\omega,a,u_1)$ is jointly continuous in $\omega, a$, the desired property holds. 
\end{Proof}

\subsection{Contraction Properties of the Reduced Equation}

We now study the one-dimensional reduced equation
\begin{equation*}0 = \mathcal{P}_s F_0(\omega,a,\psi(\omega,a)),\end{equation*} which we can rewrite as 
\begin{equation}\label{e:reduced} 0 = \mathcal{P}_s(\mathcal{L}_\omega - \mathcal{L}_{\omega_*})\left(av_*\cos(x) + \psi(a,\omega)(x)\right) +  \mathcal{P}_s(k_\omega*N(av_*\cos(x) + \psi(a,\omega)(x))),\end{equation}
since $\mathcal{P}_s(1 - \mathcal{P}) = \mathcal{P}_s$, $\mathcal{P}_s\mathcal{L}_{\omega_*} = 0$, and  $\mathcal{L}_{\omega_*}av_*\cos(x)=0$. 

We would like to find a one-parameter family of solutions $(\omega,a)$ to this equation near $(\omega_*,0)$. Typically, one would use the Implicit Function Theorem; however, the nonlinearity $N$ is only $C^1$, and both first partial derivatives of the right hand side vanish at $(\omega_*,0)$. Because $a = 0$ is a solution of \eqref{e:reduced} for any $\omega$, the entire equation can be divided by $a$, but since the nonlinearity $N$ is only $C^1$, the resulting equation is then only continuous. We thus use a direct contraction argument. 

\paragraph{Setup of Contraction Argument.}
We divide  \eqref{e:reduced} by $a$ and claim that we obtain an equation of the form 
 \begin{equation}\label{e:contractionsetup}
    \omega - \omega_* = R(\omega,a)
\end{equation} 
for some function $R(\omega,a)$. In fact, the principal term, after dividing, is $\mathcal{P}_s(\mathcal{L}_\omega - \mathcal{L}_{\omega_*})\left(v_*\cos(x)\right)$. We would like to identify the linear term in $(\omega - \omega_*)$ and show that it does not vanish. We find that the linear term in $\mathcal{P}_s(\mathcal{L}_\omega - \mathcal{L}_{\omega_*})\left(v_*\cos(x)\right)$ is $\alpha(\omega - \omega_*) =\mathcal{P}_s\frac{d}{d\omega}\left( k_\omega * A \cos(x) \right)\big|_{\omega = \omega_*} \cdot v_* (\omega - \omega_*)$. Then, provided the coefficient $\alpha$ is nonzero, equation \eqref{e:reduced} can be rearranged to the form \eqref{e:contractionsetup}. 


 

\begin{Proposition} The linear coefficient $\alpha  =\mathcal{P}_s\frac{d}{d\omega}\left( k_\omega * A \cos(x) \right)\big|_{\omega = \omega_*} \cdot v_* $ does not vanish, under the assumption that $d(i\omega_*) = 0, d'(i\omega_*) \neq 0$, and $k(-x) = k(x)$.
\end{Proposition}

\begin{Proof}
 By changing variables, one can calculate that 
 
 \begin{align*}\frac{d}{d\omega}\left( k_\omega * A \cos(x) \right)\big|_{\omega = \omega_*} \cdot v_* &= \frac{d}{d\omega}\left( \frac{1}{2}\left(Ae^{ix}\int_\R k(y)e^{-i\omega y}dy  + Ae^{-ix} \int_\R k(y)e^{i\omega y}dy \cdot \right)\right)|_{\omega = \omega_*}\cdot v_*; \\ &= \frac{d}{d\omega}\left( \frac{1}{2\sqrt{2\pi}}\left( \wh{k}(-i\omega)Ae^{ix} + \wh{k}(i\omega)Ae^{-ix} \right)\right)|_{\omega = \omega_*}\cdot v_*; \\ 
&=  \frac{d}{d\omega}\left( \frac{1}{\sqrt{2\pi}} \wh{k}(i\omega)A \right)\big|_{\omega = \omega_*} \cdot v_*\cos(x),\end{align*}

because $k$ is even. Then \begin{align*}
\alpha = \mathcal{P}_s \left( \frac{d}{d\omega}\left( \frac{\wh{k}(i\omega)}{\sqrt{2\pi}} A \right)\big|_{\omega = \omega_*} \cdot v_*\cos(x) \right) 
 &= \frac{1}{\sqrt{2\pi}}\left\langle \frac{d}{d\omega}(\wh{k}(i\omega)A)\big|_{\omega = \omega_*} v_*\cos(x), v_{ad}\cos(x)\right\rangle_{L^2([0,2\pi],\R^n)} \\
 &=\frac{1}{\sqrt{2\pi}}\left\langle \frac{d}{d\omega}\wh{k}(i\omega)A\big|_{\omega = \omega_*} v_*, v_{ad}\right\rangle_{\R^n}
 ,\end{align*}
 so it remains to show that the latter is nonzero.

 By the hypothesis, we have 
 \[
 d(i\omega_*) = \det(-\Id + \wh{k}(i\omega_*)A) \neq 0,\qquad \text{ and } \quad d'(i\omega_*) = \det(\wh{k}'(i\omega_*)A)\neq 0
 .\]
 Let $e_0$ be the first standard basis vector in $\R^n$. Because $d(i\omega_*) = 0$, there exists an invertible matrix $T$ such that $\ker(T(-\Id + \wh{k}(i\omega_*)A)T^{-1}) = e_0$; that is, 
 \[
 T\left(-\Id + \wh{k}(\nu)A\right)T^{-1} = \begin{pmatrix} b_1 (\nu - i\omega_*)  \ \Big| B_2 + \mathcal{O}((\nu - i\omega_*)) \end{pmatrix} 
 \]
 for $b_1$ a nonzero vector and $B_2$ a $n \times (n-1)$ matrix. Then by termwise expansion, 
\[
\det\left(T\left(-\Id + \wh{k}(\nu)\right)A)T^{-1}\right) = (\nu - i\omega_*) \det\begin{pmatrix} b_1  \Big|  B_2\end{pmatrix} + \mathcal{O}((\nu - i\omega_*)^2).
 \]
 By the assumption that $d'(i\omega_*) \neq 0$, we must have that $\det\begin{pmatrix} b_1 | B_2\end{pmatrix} \neq 0$. Then $b_1$ is not in the range of $B_2$ and therefore not in the range of $T(-\Id + \wh{k}(i\omega_*)A)T^{-1}$. Lastly, noticing that 
 \[
 b_1 = \frac{d}{d\nu}\left(T\left(-\Id + \wh{k}({\nu})A\right)T^{-1}\right)\Big|_{\nu = i\omega_*}\cdot e_0,
 \] 
 we get that $\wh{k}'(i\omega_*)Av_*$ is a nontrivial element of the cokernel of $(-\Id + \wh{k}(i\omega_*)A)$, since $e_0$ corresponds to $v_*$ in the original coordinates. This fact then implies that $\langle \frac{d}{d\omega}\wh{k}(i\omega)A\big|_{\omega = \omega_*} v_*, v_{ad}\rangle_{\R^n} \neq 0$, as desired. \end{Proof}
 
 Then, since $\alpha$, the coefficient of $(\omega - \omega_*)$, is nonzero, we rewrite the reduced equation \eqref{e:reduced} in the form \[(\omega - \omega_*) = R(\omega,a) = \mathcal{P}_s\wt{R}(\omega,a),\] where
 \begin{align*}
    \wt{R}(\omega, a) &= \frac{-1}{\alpha}\Bigg(\frac{1}{a}\Big[(k_\omega - k_{\omega_*})*\psi(\omega,a) + k_\omega*N(av_*\cos(x) + \psi(\omega,a))\Big]\\  & \qquad \qquad - (k_\omega - k_{\omega_*})*Av_*\cos(x) - (\omega - \omega_*) \frac{d}{d\omega}\big(k_\omega * Av_*\cos(x)\big)\Big|_{\omega = \omega_*}\Bigg) \\ 
    &:= \wt{R}_1(\omega,a) + \wt{R}_2(\omega,a) + \wt{R}_3(\omega) . 
    \end{align*}
\paragraph{Contraction Properties.}
The remainder of the section will be dedicated to showing that the function $R(\cdot,a)$ is a contraction mapping in $\omega$ on a sufficiently small neighborhood of $\omega_*$, for $a$ sufficiently small.
 
\begin{Proposition}\label{a_ep}  For any $\ep$ sufficiently small, there exists $a_*$ sufficiently small such that for any $a < a_*$, $R(\cdot, a)$ is a map from the interval $(\omega_*- \ep, \omega_* + \ep)$ into itself.\end{Proposition}

\begin{Proof} One can readily calculate $|R(\omega,a)| = |\mathcal{P}_s \wt{R}(\omega,a) |_\R \le 2 \lv \wt{R}(\omega,a) \rv_{L^\infty}$, so we investigate $ \lv \wt{R}(\omega,a) \rv_{L^\infty}$ for simplicity. 

Consider first $\wt{R}_1$. 
We have 
\[ \lv \wt{R}_1 \rv_{L^\infty} = 
\lv \frac{1}{\alpha} \frac{1}{a}(k_\omega - k_{\omega_*})*\psi(\omega,a)\rv_{L^\infty} \le |\omega - \omega_*|(\sup_\omega \lv \partial_\omega k_\omega \rv_{L^1}) \lv \frac{1}{a}\psi(\omega,a)\rv_{L^\infty}.
\] 
Since $\lv \partial_\omega k_\omega \rv_{L^1}$ is continuous in $\omega$ and hence bounded on a neighborhood of $\omega_*$, we would thus like to show that $\lv \frac{1}{a}\psi(\omega,a)\rv_{L^\infty}$ is small in a neighborhood of $(\omega_*,a)$. We can expand $\psi$ in $a$ at $a = 0$, noting that $\psi(\omega, 0) = 0$, to get $\psi(\omega,a) = a\partial_a\psi(\omega,0) + \psi_{1}(\omega,a)$, with the remainder term $\psi_1$ being jointly $C^1$ in $\omega,a$ and uniformly $o(a)$ on a neighborhood of $\omega_*$. We then note that ${\partial_a}\psi(\omega_*,0) = 0$, by the chain rule:  
\[
\frac{\partial\psi}{\partial a}(\omega_*,0) = (\partial_{u_1}F_1(\omega_*,0,\psi(\omega_*,0)))^{-1} \partial_a F_1(\omega_*,0,\psi(\omega_*,0)) = \mathcal{L}_{\omega_*}^{-1}(\mathcal{P}(N'(0)v_*\cos(x))) = 0.
\]
Therefore, $\partial_a \psi(\omega,0)$ is equal to $0$ at $\omega = \omega_*$, and continuous. We also have that $\frac{1}{a} \psi_1(\omega,a)$ is $o_a(1)$, uniformly in $\omega$ on a neighborhood of $\omega_*$, since $\psi_1(\omega,a)$ is locally uniformly $o(a)$. Therefore there exists $\ep_1$ such that for $|\omega - \omega_*| < \ep_1$, and $a$ sufficiently small, 
\[ \lv \wt{R}_1\rv_{L^\infty} \le \frac{1}{\alpha}|\omega - \omega_*|(\sup_\omega \lv \partial_\omega k_\omega \rv_{L^1}) \lv \frac{1}{a}\psi(\omega,a)\rv_{L^\infty} < \frac{\ep}{6}. 
\]
 As for $\wt{R}_2$, we note that 
 \[
 \lv k_\omega*N(av_*\cos(x) + \psi(\omega,a))\rv_{L^\infty} \le \lv k_\omega \rv_{L^1} \lv N(av_*\cos(x) + \psi(\omega,a))\rv_{L^\infty}.
 \]
 We know that $\lv k_\omega \rv_{L^1}$ can be bounded on a neighborhood of $\omega_*$ since it is continuous in $\omega$. Now, we have $N(0) = 0$, and $\frac{\partial}{\partial a} (N(av_*\cos(x) + \psi(\omega,a)))|_{a = 0} = N'(0)\frac{\partial}{\partial_a}\psi(\omega,a)|_{a=0} = 0$, since $N'(0) = 0$. Thus we will have $\lv N(av_*\cos(x) + \psi(\omega,a))\rv_{L^\infty} = o(a)$, uniformly in a neighborhood of $\omega_*$, since $N$ and $\psi$ are $C^1$.
 Then, given any 
 $\ep > 0$, for $a$ sufficiently small, 
 \[
 \lv \wt{R}_2 \rv_{L^\infty} \le \frac{1}{\alpha} \sup_\omega \lv k_\omega \rv_{L^1} \lv \frac{1}{a}N(av_*\cos(x) + \psi(\omega,a))\rv_{L^\infty} < \frac{\ep}{6}.
 \]
 We claim that $\wt{R}_3$ is at least quadratic in $(\omega - \omega_*)$. Since $(k_\omega*Av_*\cos(x))(\cdot) = (k*Av_*\cos(\omega x))(\frac{1}{\omega} \cdot)$, and the latter is smooth in $\omega$, we can expand in $\omega$ and find that the expansion starts at quadratic order.  Then there exists $\ep_2$ such that for $|\omega - \omega_*| < \ep_2,$ $\lv\wt{R}_3\rv_{L^\infty}$ is less than $\frac{\ep_2}{6}$.
 
 Thus, for any $\ep < \min(\ep_1,\ep_2),$ with $a$ sufficiently small, for $|\omega - \omega_*| < \ep $, \[ 
 \lv \wt{R}(\omega,a)\rv_{L^\infty} \le \lv \wt{R}_1 \rv_{L^\infty} + \lv \wt{R}_2 \rv_{L^\infty} + \lv \wt{R}_3 \rv_{L^\infty} < \frac{\ep}{6} + \frac{\ep}{6} + \frac{\ep}{6} = \frac{\ep}{2},
 \]
 $\text{so that }| R((\omega - \omega_*) + \omega_*,a)|< \ep$.
  \end{Proof}
  It remains to show that $R(\cdot,a)$ is a contraction mapping. 
 
\begin{Lemma}\label{contraction}
  There exists $\ep >0$ and $a_*>0$ such that for $a < a_*$, the map $R(\omega,a)$ is a contraction mapping from $(\omega_* -\ep, \omega_* + \ep)$ to itself. 
 \end{Lemma}
 \begin{Proof}
 We investigate the Lipschitz constant of $R$. We have 
 \begin{align*}
    R(\omega_1, a) - R(\omega_2,a) = \frac{-1}{\alpha}\frac{1}{a}\mathcal{P}_s&\Bigg[(k_{\omega_1} - k_{\omega_2})*\psi(\omega_2,a) + (k_{\omega_1} - k_{\omega_*})*(\psi(\omega_1,a) - \psi(\omega_2,a)) \\ &+ (k_{\omega_1} - k_{\omega_2})*N(av_*\cos(x) + \psi(\omega_2,a)) \\ &+k_{\omega_1}*\left(N(av_*\cos(x) + \psi(\omega_1,a)) - N(av_*\cos(x) + \psi(\omega_2,a))\right)\\  &+ (k_{\omega_1} - k_{\omega_*})*Av_*\cos(x) - (\omega_1 - \omega_*) \frac{d}{d\omega}\left(k_\omega * Av_*\cos(x)\right)\big|_{\omega = \omega_*} \\ 
    &- \left((k_{\omega_2} - k_{\omega_*})*Av_*\cos(x) - (\omega_2 - \omega_*) \frac{d}{d\omega}\left(k_\omega * Av_*\cos(x)\right)\big|_{\omega = \omega_*}\right)\Bigg]. 
\end{align*}
We again estimate norms in $L^\infty$, accounting for the factor of 2. For the first term, as before, we have 
\[
    \lv \frac{1}{a} (k_{\omega_1} - k_{\omega_2})*\psi(\omega_2,a)\rv_{L^\infty} \le |\omega_1 - \omega_2|\left(\sup_\omega \lv \partial_\omega k_\omega \rv_{L^1}\right) \lv \frac{1}{a}\psi(\omega_2,a)\rv_{L^\infty}
.\] 
By the argument above, for the values of $(\omega,a)$ considered, we already have 
\[
\frac{1}{\alpha}(\sup_\omega \lv \partial_\omega k_\omega \rv_{L^1}) \lv \frac{1}{a}\psi(\omega,a)\rv_{L^\infty} < \frac{1}{6}, \qquad \text{ so that}\quad  \lv \frac{1}{\alpha} \frac{1}{a} (k_{\omega_1} - k_{\omega_2})*\psi(\omega_2,a)\rv_{L^\infty} \le \frac{1}{6}|\omega_1 - \omega_2|.
\]
    
For the second term, we have \[\lv \frac{1}{\alpha} (k_{\omega_1} - k_{\omega_*})*(\psi(\omega_1,a) - \psi(\omega_2,a))\rv_{L^\infty} \le |\omega_1 - \omega_2| \cdot |\omega_1 - \omega_*|\frac{1}{\alpha} (\sup_\omega \lv \partial_\omega k_\omega \rv_{L^1})\frac{1}{a}\lip_\omega(\psi(\omega,a)).\] The term $\lip_\omega(\psi(\omega,a))$ is bounded by $|\partial_\omega\psi(\omega,a)|$, which is bounded by $C_1|a|$, so that $\frac{1}{a} \lip_\omega(\psi(\omega,a)) \le C_1$. Then for $|\omega_1 - \omega_*|$ sufficiently small, the whole term will have small Lipschitz constant: there exists $\ep_3$ such that for $|\omega_1 - \omega_*| < \ep_3$, 
\[
\lv \frac{1}{\alpha} (k_{\omega_1} - k_{\omega_*})*(\psi(\omega_1,a) - \psi(\omega_2,a))\rv_{L^\infty} < \frac{1}{12}|\omega_1 - \omega_2|.
\]
For the third term, we have 
\begin{align*}
\lv \frac{1}{\alpha} (k_{\omega_1} - k_{\omega_2})*N(av_*\cos(x) + \psi(\omega_2,a)) \rv_{L^\infty} &\le |\omega_2 - \omega_1|\frac{1}{\alpha}\sup_\omega \frac{d}{d\omega} \lv k_\omega *N(av_*\cos(x) + \psi(\omega_2,a)) \rv_{L^\infty}\\
&\le |\omega_1 - \omega_2|\frac{1}{\alpha}\sup_\omega \lv \partial_\omega k_\omega \rv_{L^1}\lv N(av_*\cos(x) + \psi(\omega_2,a))\rv_{L^\infty} 
.\end{align*}
As discussed previously, $\frac{1}{a}\lv N(av_*\cos(x) + \psi(\omega_2,a))\rv_{L^\infty} $ is $o_a(1)$, and the rest of the terms are bounded. Then for $a$ sufficiently small, $\lv \frac{1}{\alpha} (k_{\omega_1} - k_{\omega_2})*N(av_*\cos(x) + \psi(\omega_2,a)) \rv_{L^\infty} < \frac{1}{12}|\omega_1 - \omega_2|$. 
 
 For the fourth term, note that \begin{align*} \lip_\omega (\frac{1}{a}N(av_*\cos(x) + \psi(a,\omega))) &\le \frac{1}{a}\sup |N^{'}(av_*\cos(x) + \psi(\omega,a))| \lip_\omega \psi  \\ 
 &\le \frac{1}{a}\sup |N^{'}(av_*\cos(x) + \psi(\omega,a))| C_1 |a| \\
 &\le C_1\sup |N^{'}(av_*\cos(x) + \psi(\omega,a))|.\end{align*}
 
 Because $N'(0) = 0$, with $N'$ continuous, and the argument  $av_*\cos(x) + \psi(\omega,a)$ equals $0$ at $(\omega_*,0)$, there exists a neighborhood of $(\omega_*, 0)$ for which $\frac{1}{\alpha} C_1 \lv k_{\omega_1} \rv_{L^1} \sup |N^{'}(av_*\cos(x) + \psi(\omega,a))| < \frac{1}{12}$. Then on that neighborhood, 
 \[
 k_{\omega_1}*\left(N(av_*\cos(x) + \psi(\omega_1,a)) - N(av_*\cos(x) + \psi(\omega_2,a))\right) < \frac{1}{12}|\omega_1 - \omega_2|.
 \]
  As to the last difference of terms, which is independent of $a$, note that because the term is quadratic in $(\omega - \omega_*)$, then there exists $\ep_4$ for which the Lipschitz constant is less than $\frac{1}{12}$ for $|\omega - \omega_*| < \ep_4$. 
 
 Now, fix $\ep < \min(\ep_1, \ep_2, \ep_3,\ep_4)$, and $a_0$ sufficiently small such that $(-a_0, a_0) \times (-\ep, \ep)$ is in all neighborhoods mentioned above, and so that for $a < a_0$, by Lemma \ref{a_ep}, the map $R(\cdot, a)$ maps $B_\ep(\omega_*)$ to itself. Then we can find $a_*$, possibly smaller, such that for $a < a_*$, 
 \[\lv \wt{R}(\omega_1, a) - \wt{R}(\omega_2, a) \rv_{L^\infty} < (\frac{\ep}{6} + \frac{\ep}{12} + \frac{\ep}{12}+ \frac{\ep}{12}+ \frac{\ep}{12})|\omega_1 - \omega_2| = \frac{\ep}{2} |\omega_1 - \omega_2|, 
 \] 
 so that 
 \[
 |R(\omega_1, a) - R(\omega_2,a)| < \ep|\omega_1 - \omega_2|.
 \] 
 Then for all $a < a_*$, $R(\cdot,a)$ is a contraction in $(\omega - \omega_*)$ on $(-\ep, \ep)$. 
 \end{Proof} 
 
 \subsection{Proof of Theorems \ref{per_existence} and \ref{LC}}

Using the above contraction properties, we can now prove Theorems \ref{per_existence} and \ref{LC}. 

\begin{Proof}[(of Theorem \ref{per_existence})]
We show existence of a two-parameter family of solutions to equation \eqref{e:nonlin}. 
 
Let $a < a_*$. Then by Lemma \ref{contraction} and the Banach fixed point theorem, there exists a unique fixed point of $\omega = R(\omega,a)$. Then for all $a < a_*$, there exists $\omega(a)$ such that $F(\omega(a), a, \psi(\omega(a),a)) = 0$.
  
The family of solutions $u(x;a) = av_*\cos(x) + \psi(\omega(a),a)(x)$ is then a one-parameter family of solutions to \eqref{e:nonlinper} near $a = 0$, which in turn yields a family of solutions $u_c(x;a) = u(\omega(a) x;a)$ to the original equation \eqref{e:nonlin}. In order to obtain a two-parameter family of solutions, we use the fact that the original equation \eqref{e:nonlin} is translation-invariant. This ensures that the function $u_c(\cdot + \tau;a)$ is a solution for any $\tau$. Lastly, the properties $u_c(y;0) = 0, \omega(0) = \omega_*$ are easily verified by examining properties of $\psi$ and $R$. \end{Proof}

This establishes a two-parameter set of periodic solutions to \eqref{e:nonlin}. However, to prove Theorem \ref{LC}, we need to further characterize this set of solutions topologically: 

\begin{Proposition}\label{2dparam}
 There exists a neighborhood of the origin $U_p\subset \R^2$ and a continuous map
 $S:U_p\to C^0_{-\eta}(\R,\R^n)$ whose range consists of the family of continuous periodic solutions to \eqref{e:nonlin} found in Theorem \ref{per_existence}. 
 \end{Proposition}
\begin{Proof} 
First, identifying $\R^2$ with $\C$, let $s_1: \C \backslash \{0\} \to \R_{>0} \times [0,2\pi)$ be defined by $s_1(z) = (|z|, \mathrm{arg}(z))$, and let $s_2: \R_{>0} \times [0,2\pi) \to C^0_{-\eta}(\R,\R^n)$ be defined by $s_2(r, \theta) = u_c(x + \theta;r)$. The composition $s_2 \circ s_1$ can be seen to be continuous and one-to-one from $\C \backslash \{ 0 \}$ to $C^0_{-\eta}(\R,\R^n)$. Then we would like to extend $s_2 \circ s_1$ continuously to $\C$. Let 
\[
S(z) = \begin{cases} (s_2 \circ s_1)(z), & z \neq 0, \\ 0, & z = 0. \end{cases}
\] 
We can see that $S$ is continuous at 0 because $\lv u_c(\cdot + \arg(z);|z|)\rv_{C^0}$ approaches 0 as $|z|$ approaches 0; $s$ is also still one-to-one. Since there exists $a_* > 0 $ such that $u_c(\cdot + \omega(a)\tau,a)$ is a solution to \eqref{e:nonlinper} for any $a > 0, \tau \in \R$, then  $S$ is a continuous, one-to-one map from the neighborhood $\{|z| < a_*\}$ to the set of continuous periodic solutions to \eqref{e:nonlin}. 
\end{Proof}

With this characterization, we now turn to the proof of Theorem \ref{LC}.
\begin{Proof}[(of Theorem \ref{LC})] Letting $\wt{N}(u) = A^{-1}N(u)$, $K(x) = k(x)\cdot A$, then equation \eqref{e:nonlin} is in the appropriate form for Theorem \ref{CMexistence}. The kernel $\mathcal{E}_0$ of $\T u = -u + K*u$ in this case is two-dimensional, since $d(\nu)$ has two single roots on the imaginary axis. Then by Theorem \ref{CMexistence}, there exists $\delta > 0$, a center manifold $\mathcal{M} \subset C^0_{-\delta}$ and a map $\Psi: \mathcal{E}_0 \to \mathcal{M}$, with $\Psi(0) = 0, D\Psi(0) = 0$, such that $\mathcal{M} = \{ u_0 + \Psi(u_0)\ | \ u_0 \in \ker\T\}$. By property $(iv)$ of the theorem, $\mathcal{M}$ contains all solutions $u$ to \eqref{e:nonlin} with $\lv u \rv_{C^0} < \ep$ for some $\ep > 0$. Then, taking $\ep_* < \min(\ep, a_*)$, the family of solutions $\{u_c(\cdot + \tau;a) \ |  \ \tau \in \R, a \in [0,\ep_*) \}$ is contained, as a set, in $\mathcal{M}$. 

The composition of maps $\mathcal{Q} \circ S$, where $\mathcal{Q}$ is the projection onto $\mathcal{E}_0$ as defined in Section 4, and $S$ is the map from Proposition \ref{2dparam}, is then a continuous, one-to-one map from the neighborhood $U_p$ of 0 in $\R^2$ to $\mathcal{E}_0$. Note that $\mathcal{Q}$ is one-to-one because it is invertible on $\mathcal{M}$. Its restriction to a closed neighborhood of 0 contained in $U_p$ will therefore have continuous inverse and hence be open. Then the image of $U_p$ in $\mathcal{E}_0$ contains a ball of positive radius in $\mathcal{E}_0$, which, since $\mathcal{E}_0$ is finite-dimensional, contains a ball in $\mathcal{E}_0$ under the $C^0$ norm. Lastly, since $\lv \mathcal{Q}u\rv_{C^0} \le \lv u \rv_{C^0}$, any solution to \eqref{e:nonlin} with sufficiently small $C^0$ norm is in the image of $U_p$ in the $\mathcal{M}$. Hence any sufficiently small solution to \eqref{e:nonlin} is periodic, which proves Theorem \ref{LC}. \end{Proof}

\section{Discussion}

We have established Fredholm properties for a nonlocal operator with a multiplication operator as its principal part, finding an additional source of noncompactness corresponding to zeros of the principal part. Using this theory, we established existence of finite-dimensional center manifolds for nonlocal equations on $C^0$-based spaces, allowing for optimal regularity of the manifold in a set of coordinates. This allowed us to prove a nonlocal Lyapunov-Center theorem in the $C^1$ case. We describe briefly below possible further directions of this work, and some apparent difficulties therein. 

\paragraph{General Nonlinearities.} The work here focuses on pointwise substitution operators as a simple class allowing for optimal regularity; a natural extension is to consider general Frechet operators on function spaces. One limitation is establishing the bootstrapping step for these operators, which involves smoothness of the inverse of $(\mathrm{Id} + \mathcal{G}^\ep)$. 

\paragraph{Optimal regularity without changing coordinates.} A natural question is whether optimal regularity can be obtained in the original equation without changing variables, possibly in different function spaces. The inherent difficulty is that differentiating the shift operator requires that the trajectory be differentiable. It is not clear how regularity could be obtained in these coordinates using for instance bootstrapping. On the other hand, it seems plausible that vector fields are simply optimally regular only in this particular choice of coordinates: changing coordinates for an ODE with $C^1$ vector field with a $C^1$ diffeomorphism of course only results in a continuous vector field, albeit with a well defined $C^1$ flow. 


\paragraph{Extension to a Cylinder.} The systems studied here are in one spatial variable $\xi \in \R$. As in local spatial dynamics, one would like to extend the theory to the 2-dimensional, cylindrical case (as in \cite{kirchgassner1982wave} by Kirchgassner). One would have to find conditions under which the kernel of the linearization is finite-dimensional. Much loftier and less clear, but no less interesting, would be an extension to 2 or more unbounded spatial variables, where the time-like flow would correspond to a more general symmetry.

\paragraph{Regularity of the Kernel.} The present argument relies on regularity of the convolution kernel---enough to map $L^p$ into $W^{1,p}$. It is conceivable that this assumption could to be relaxed slightly, such as to a kernel mapping $L^p$ to $W^{\theta,p}, \theta >0$, exploiting repeated bootstrapping. 

\paragraph{Localization of the Kernel.} Computing Fredholm indices and constructing center manifolds requires exponential localization of the kernel. Inspecting however the way multiplicities and crossing numbers are computed, or the way Taylor expansions of reduced vector fields are determined, one finds that only finite, possibly high moments of the kernel enter the calculation. One may therefore suspect that moment conditions would be sufficient to establish some, possibly weaker result. It seems however difficult to guarantee the robustness with respect to parameters and the fact that center manifolds contain all bounded solutions without such strong localization assumptions (or additional structure such as monotonicity). Existence of small bounded solutions alone, can indeed be deduced from appropriate moment conditions alone in many scenarios; see for instance \cite{ST19}



\paragraph{Extension to Other Function Spaces.} Lastly, the choice of $C^0$-based spaces here was a natural choice of spaces where pointwise nonlinearities do not lose regularity as substitution operators. Regularity questions when studying for instance equations in cylindrical domains may well require different function spaces, such as spaces with H\"older regularity. It is conceivable that the strategy pursued here may well generalize, although cut-off procedures may be more involved. 


\section*{Declarations}
 
\paragraph{Ethical Approval.}
This declaration is not applicable.
\vspace{-0.2 in}
\paragraph{Competing interests.}
The authors declare that they have no competing interests as defined by Springer, or other interests that might be perceived to influence the results and/or discussion reported in this paper.
\vspace{-0.2 in}
\paragraph{Authors' contributions.}
All authors contributed equally to the manuscript.
\vspace{-0.2 in}
\paragraph{Funding.}
The authors acknowledge partial support through grants NSF DMS-1907391 and NSF DMS-2205663.
\vspace{-0.2 in}
\paragraph{Availability of data and materials.}
This declaration is not applicable.

\paragraph{Acknowledgments.} We thank the anonymous referee for careful reading and many constructive suggestions to improve the manuscript.



\end{document}